\documentclass[11pt]{article}
\usepackage{eurosym}
\usepackage{geometry}
\usepackage{graphicx}
\usepackage{amssymb,amsmath}
\usepackage{theorem}
\usepackage{epstopdf}

\newtheorem{theo}{Theorem}
\newtheorem{lem}{Lemma}

\newtheorem{rem}{Remark}
\numberwithin{equation}{section}
\numberwithin{lem}{section}
\numberwithin{prop}{section}
\numberwithin{theo}{section}
\numberwithin{rem}{section}
\begin{document}

\title{Strong Gaussian approximation of the mixture Rasch model}
\author{Friedrich Liese%
\thanks{%
Institut f\"ur Mathematik, Universit\"at Rostock, D-18051 Rostock, Germany,
\textsc{email}: friedrich.liese@uni-rostock.de} \hspace{2.5cm} Alexander
Meister\thanks{%
Institut f\"ur Mathematik, Universit\"at Rostock, D-18051 Rostock, Germany,
\textsc{email}: alexander.meister@uni-rostock.de} \hspace{2.5cm}
Johanna Kappus\thanks{%
Institut f\"ur Mathematik, Universit\"at Rostock, D-18051 Rostock, Germany,
\textsc{email}: johanna.kappus@uni-rostock.de}}
\maketitle

\begin{abstract}
We consider the famous Rasch model, which is applied to psychometric surveys when $n$ persons under test answer $m$ questions. The score is given by a realization of a random binary $n\times m$-matrix. Its $(j,k)$th component indicates whether or not the answer of the $j$th person to the $k$th question is correct. In the mixture Rasch model one assumes that the persons are chosen randomly from a population. We prove that the mixture Rasch model is asymptotically equivalent to a Gaussian observation scheme in Le Cam's sense as $n$ tends to infinity and $m$ is allowed to increase slowly in $n$. For that purpose we show a general result on strong Gaussian approximation of the sum of independent high-dimensional binary random vectors. As a first application we construct an asymptotic confidence region for the difficulty parameters of the questions.  
\end{abstract}



\section{Introduction}


The Rasch model is a famous and widely used approach to analyse surveys in
the field of psychometrics. It assumes that each of $n$ subjects (typically
persons) are exposed to $m$ items (typically questions to be answered). For
each $j=1,\ldots,n$ and $k=1,\ldots,m$ the correctness of the answer of person $j$ to the question $k$ is a binary random variable $X_{j,k}$ where the probability of
a correct answer, i.e. $X_{j,k}=1,$ is given by%
\begin{equation*}
P(X_{j,k}=1)\,=\,\frac{\exp \{\beta_{j}-\theta_{k}\}}{1+\exp \{\beta_{j}-\theta_{k}\}}\,,\qquad k=1,\ldots ,m;\,j=1,\ldots ,n\,.
\end{equation*}%
The parameter $\theta _{k}$ characterizes the difficulty of the $k$th item
and the parameter $\beta _{j}$ reflects the ability of the $j$th
individual. The $\beta _{j}$ may be either considered as unknown parameters (standard Rasch model) 
or as realizations of i.i.d. random variables with distribution $F$. The
latter case describes the situation in which the individuals are randomly selected
from a large population. Then the observation vectors $Y_{j}=(X_{j,1},\ldots,X_{j,m})$
are i.i.d. and it holds for every binary matrix $\varepsilon =(\varepsilon
_{j,k})_{j=1,\ldots ,n;k=1,\ldots ,m\,}$ that
\begin{equation}
P\left( (Y_{1},\ldots,Y_{n})^T = \varepsilon \right) =\prod\nolimits_{j=1}^{n}\int
\left\{ \prod\nolimits_{k=1}^{m}\frac{\exp \{\varepsilon _{j,k}(\beta
-\theta _{k})\}}{1+\exp \{\beta -\theta _{k}\}}\right\} dF(\beta )\,,\qquad
\varepsilon _{j,k}\in \{0,1\}.  \label{eq:0.1}
\end{equation}%
This type of psychometric model is called the mixture Rasch model 
which will be the central object in this paper.

For original literature we refer to the book of Rasch (1960/1980), after
whom the model has been named. Also we mention the books of Alagumalai et
al. (2005) and Bezruczko (2005) for applications of the Rasch model. It has also confined attention in the econometric literature
(Hoderlein et al. (2011)). The mixture model is used in Lindsay et al. (1991), Rice (2004) and Strasser
(2012a,b). Also we refer to the books of Fischer and Molenaar (1995) and
von Davier and Carstensen (2007).

So far most of the literature on the Rasch model has mainly
focused on the estimation of the difficulty parameters, consistency and
asymptotic normality for bounded $m$ where maximum likelihood (ML) or
quasi-ML methods are preferred, see e.g. de Leeuw and Verhelst (1986) or
Pfanzagl (1993, 1994). Lindsay et al. (1991) consider semiparametric
estimation in the Rasch model and related problems. Biehler et al. (2015)
study saddlepoint approximation of the ability parameters. Doebler et al.
(2013) construct confidence intervals for the ability parameters. Strasser
(2012a,b) thoroughly investigates the covariance structure and asymptotic
distribution of quasi-ML estimators in the mixture Rasch model.

In this work we approximate the mixture
Rasch model in the strong Le Cam sense by a model which contains a Gaussian
observation, and -- conditionally on that -- another Gaussian observation  whose distribution does not depend on the ability distribution $F$ (as $n\rightarrow \infty $). This investigation is motivated by the fact that, for
Gaussian models, the structure of optimal estimators and tests is understood very well in both the parametric and nonparametric case.

As a first application we will construct a uniform asymptotic confidence
ellipsoid for the difficulty parameters in the asymptotically equivalent
Gaussian model under potentially increasing (but restricted) dimension $m$,
which, thus, also represents a uniform asymptotic confidence ellipsoid in
the original mixture Rasch model. Also the asymptotic equivalence result
will open a broad field of further applications as we will explain in the
conclusions.

The distribution $F$ in (\ref{eq:0.1}) is not nonparametrically identified for bounded $m$, a situation that is similar to the binomial mixture models.
Therefore we allow $m=m_{n}$ to tend to infinity, as $n\rightarrow \infty $. Therein $m_{n}$ has to be of smaller order compared to $n$. This means that \ there are much more subjects under test compared
to the total number of questions contained in the sheet, a condition that is
satisfied in almost all applications and especially in the Programme for
International Student Assessment (PISA), to which the Rasch model has been
applied.

\section{Asymptotic Equivalence} \label{aseq}

In this section we provide a brief introduction to the concept of asymptotic
equivalence. Assume we have two statistical experiments $\mathcal{E}%
_{j}=\left( \Omega _{j},\mathfrak{A}_{j},(P_{\theta ,j})_{\theta \in \Theta
}\right) $ with the same parameter space $\Theta $. By $\mathfrak{K}_{i,j}$
and $i,j\in \{1,2\}$ we shall denote the set of all Markov  kernels $%
K_{i,j}:\mathfrak{A}_{j}\times \Omega_{i}\rightarrow \lbrack 0,1]$. The
application of $K_{2,1}$ on $P_{\theta,2},$ i.e.
\begin{equation*}
(K_{2,1}P_{\theta,2})(A_{1})=\int K_{2,1}(A_{1},\omega _{2})P_{\theta
,2}(d\omega _{2})
\end{equation*}%
is a probability measure on $\left( \Omega_{1},\mathfrak{A}_{1}\right)$. 
The two statistical experiments are called equivalent if there are Markov kernels $K_{1,2}$ and $K_{2,1}$, both not depending on $\theta$, such that $K_{2,1}P_{\theta,2}=P_{\theta,1}$ and $K_{1,2}P_{\theta,1}=P_{\theta,2}$ for all $\theta \in \Theta$.  Then the two experiments
are also equivalent in the decision theoretic sense. Indeed, if $(\mathcal{D},\mathfrak{D})$ is a decision space, $L(a,\theta)$ a loss function and $D_{i}(B,\omega _{i}),B\in \mathfrak{D},\omega _{i}\in \Omega_{i}$ is a (randomized) decision for the $i$th  experiment then
\begin{equation*}
D_{j}(B,\omega_{j}):=\int D_{i}(B,\omega_{i})K_{j,i}(d\omega_{i},\omega_{j})
\end{equation*}%
is a decision for the other experiment and it can be easily seen that both
decisions have identical risk functions. Now suppose that $T:\Omega_{1}\rightarrow \Omega_{2}$ is sufficient, i.e. there exists some Markov kernel $K$ which does not depend on $\theta$ but represents a version of the conditional measure given $T$ under $P_{\theta,1}$ for all $\theta\in\Theta$; concretely $K(A_1,T) = E_{\theta,1}(1_{A_1} \mid T)$ $P_{\theta,1}$-a.s., for all $A_1 \in \mathfrak{A}_1$; and that $P_{\theta
,2} = P_{\theta ,1}\circ T^{-1}$ for all $\theta\in\Theta$. Let $\delta_{a}$ denote the Dirac measure
concentrated at point $a.$ Then $K_{1,2}(A_{2},\omega _{1})=\delta
_{T(\omega _{1})}(A_{2})$ is a Markov kernel and it holds that $%
K_{1,2}P_{\theta ,1}=P_{\theta ,1}\circ T^{-1}.$ The sufficiency of $T$
implies that there is a Markov  kernel $%
K_{2,1}$ with $K_{2,1}P_{\theta,2}=P_{\theta,1}.$ The two models are
equivalent, therefore.

The concept of deficiency makes precise in what sense the approximate
sufficiency of a statistic or, more generally, the approximate equivalence is to
be understood. It is defined with the help of the total variation distance $\mathsf{TV}(P,Q)=2 \sup_{A}|P(A)-Q(A)|$ between the distributions $P$ and $Q$. Put, for $i,j\in
\{1,2\},i\neq j$, 
\begin{eqnarray*}
\delta (\mathcal{E}_{i},\mathcal{E}_{j}) &=&\inf_{K_{j,i}\in \mathfrak{K}%
_{i,j}}\sup_{\theta \in \Theta }\mathsf{TV}(K_{j,i}P_{\theta ,j},P_{\theta
,i}) \\
\Delta (\mathcal{E}_{1},\mathcal{E}_{2}) &=&\max (\delta (\mathcal{E}_{1},%
\mathcal{E}_{2}),\delta (\mathcal{E}_{2},\mathcal{E}_{1})).
\end{eqnarray*}%
Therein $\delta (\mathcal{E}_{i},\mathcal{E}_{j})$ is called the deficiency of $%
\mathcal{E}_{i}$ and $\mathcal{E}_{j}$ and $\Delta (\mathcal{E}_{1},\mathcal{%
E}_{2})$ is the Le Cam distance of $\mathcal{E}_{1}$ and $\mathcal{E}_{2}.$
It is a metric in the space of equivalence classes of statistical
experiments with a joint parameter set. Two sequences $\mathcal{E}%
_{j,n}=\left( \Omega _{j,n},\mathfrak{A}_{j,n},(P_{\theta ,j,n})_{\theta \in
\Theta _{n}}\right) ,j=1,2$ of statistical experiments are called   
asymptotically equivalent if $\lim_{n\rightarrow \infty }\Delta (\mathcal{E}%
_{1,n},\mathcal{E}_{2,n})=0$. By a slight abuse of language one calls the experiments $\mathcal{E}_{1,n}$ and $\mathcal{E}_{2,n}$ asymptotically equivalent while this means asymptotic equivalence of the corresponding sequences. Sometimes the sample spaces are identical, then
\begin{equation*}
\Delta (\mathcal{E}_{1,n},\mathcal{E}_{2,n})\leq \sup_{\theta \in \Theta
_{n}}\mathsf{TV}(P_{\theta ,1,n},P_{\theta ,2,n}).
\end{equation*}%
Asymptotic equivalence allows to take over asymptotic properties such as
convergence rates of estimators or asymptotic confidence regions from one
experiment to the other.

In the local asymptotic decision theory $\Theta $ is an open subset of $%
\mathbb{R}^{d}$ and for a fixed $\theta _{0}\in \Theta $ and a sequence $%
a_{n}$ tending to zero one introduce a local parameter $h\in \mathcal{H}%
_{n}=\{h:\theta _{0}+a_{n}h\in \Theta \} \subseteq \mathbb{R}^{d}$. The so
called LAN\ condition for $\mathcal{E}_{1,n}$, see Strasser (1985) is
equivalent to the following statement: There is a matrix $\mathcal{I}(\theta
_{0}),$ called information matrix, such that $(\mathcal{E}_{1,n})_n$ converges
weakly to the Gaussian experiment $\mathcal{E}_{2}=(\mathbb{R}^{d},\mathfrak{B}(\mathbb{R}^{d}),(\mathsf{N}(\mathcal{I}(\theta _{0})h,\mathcal{I}(\theta _{0}))_{h\in \mathbb{R}^{d}})$. Weak convergence means that  $\Delta (\mathcal{E}_{1,n}^{%
\mathcal{H}},\mathcal{E}_{2}^{\mathcal{H}})\rightarrow 0,$ where the
superscript $\mathcal{H}$ means that for consider the experiments only for a
finite but arbritrary subset $\mathcal{H}\subseteq \mathcal{H}_{n}$ as
parameter set. A typical situation, in which this condition holds, occurs if
the family $P_{\theta }$ is $L_{2}$-differentiable and $a_{n}=\frac{1}{\sqrt{%
n}}$ and $P_{\theta ,1,n}=P_{\theta ,1}^{\otimes n},$ i.e. if we have i.i.d. observations.

For books on Le Cam theory we refer to Le Cam (1986), Strasser (1985) and Le
Cam and Yang (2000), Shiryaev and Spokoiny (2000), Liese and Miescke (2008). In
nonparametric literature research mainly focuses on showing asymptotic
equivalence of curve estimation problems to white noise models, in which the
target curve occurs as the drift function of a Wiener process. Therein we
mention e.g. Nussbaum (1996) and Carter (2002) for density estimation; Brown
and Low (1996), Rohde (2004), Carter (2006), Cai and Zhou (2014) and
Schmidt-Hieber (2014) for nonparametric regression; Meister (2011) for
functional linear regression; Rei{\ss } (2011), Genon-Catalot and Lar\'{e}do
(2014) and Mariucci (2016) for the analysis of more complex stochastic
processes. The paper of Meister and Rei{\ss } (2013) somehow deviates from
this list as it establishes asymptotoic equivalence of nonregular
nonparametric regression and a specific Poisson point process. Still
Gaussian limit models are most popular.

\section{Dimension Reduction}


The sample space for the Rasch model is $\Omega =\{0,1\}^{n\times m},$ the
space of all binary $n\times m$-matrices $\omega =(\omega_{j,k}),1\leq
j\leq n;1\leq k\leq m$. Throughout we equip a discrete sample
space $\mathcal{X}$ by the power set ${\cal P}(\mathcal{X})$ as $\sigma$-algebra. Therein $X_{j,k}(\omega )=\omega _{j,k}$ indicates the correctness of the answer of person
$j$ to question $k$. Then $Y_{j}=(X_{j,1},\ldots,X_{j,m})$ is the response vector of person $j$.

We fix some $R>0$ and a set ${\cal F}$ of admitted distributions $F$ in (\ref{eq:0.1}) and set 
\begin{equation} \label{eq:Theta} \Theta = \Big\{\theta\in [-R,R]^{m} : \sum_{k=1}^m \theta_k = 0\Big\}\,. \end{equation}
Note that the condition that the $\theta_k$ add to zero is a common calibration to ensure identifiability of the difficulty parameters. For $\theta =(\theta
_{1},\ldots,\theta _{m})$ and $F\in {\cal F}$ we denote by $P_{\theta ,F}^{A}$
the joint distribution of $Y_{1},\ldots,Y_{n}.$ The density $dP_{\theta
,F}^{A}/d\kappa _{n\times m}$ of $P_{\theta ,F}^{A}$ with respect to the
counting measure $\kappa _{n\times m}$ on $\Omega $ is the probability mass
function and (\ref{eq:0.1}) yields
\begin{equation}
\frac{dP_{\theta ,F}^{A}}{d\kappa _{n\times m}}(\omega )=P_{\theta
,F}^{A}(\{\omega \})=\prod\nolimits_{j=1}^{n}\int \left\{
\prod\nolimits_{k=1}^{m}\frac{\exp \{\omega _{j,k}(\beta -\theta _{k})\}}{%
1+\exp \{\beta -\theta _{k}\}}\right\} dF(\beta ).  \label{A}
\end{equation}%
Putting together all components we arrive at the experiment (or mixture
Rasch model)
\begin{equation*}
\mathcal{A}_{n,m}:=\left( \{0,1\}^{n\times m},{\cal P}(\{0,1\}^{n\times
m}),(P_{\theta ,F}^{A})_{(\theta ,F)\in \Theta \times {\cal F}
}\right).
\end{equation*}%
We set 
\begin{equation} \label{eq:neu}
S_{j}=\sum_{k=1}^{m} X_{j,k}, \qquad N_{k}=\sum\nolimits_{j=1}^{n} 1_{\{k%
\}}(S_{j}), \qquad T_{k}=\sum_{j=1}^{n} X_{j,k},
\end{equation}
where $1_{A}$ is the indicator function of the set $A$; and
\begin{equation*}
G(k,\theta ,F)=\log \Big\{ \int \Big\{ \exp \{k\beta
\}\prod\nolimits_{l=1}^{m}\frac{1}{1+\exp \{\beta -\theta _{l}\}}\Big\}
dF(\beta )\Big\}
\end{equation*}%
for $k=0,\ldots ,m$. The representation (\ref{A}) yields
\begin{equation*}
\frac{dP_{\theta ,F}^{A}}{d\kappa _{n\times m}}=\exp \Big\{ -
\sum\nolimits_{k=1}^{m}\theta
_{k}T_{k}+\sum\nolimits_{j=1}^{n}G(S_{j},\theta ,F)\Big\} = \exp \Big\{ - 
\sum\nolimits_{k=1}^{m}\theta
_{k}T_{k}+\sum\nolimits_{k=0}^{m}N_{k}\cdot G(k,\theta ,F)\Big\} .
\end{equation*}%
Then, by the Fisher-Neyman factorization\ criterion\ in standard Polish experiments, we realize in a first
step that the statistic $(S_{1},\ldots,S_{n},T_{1},\ldots,T_{m})$ which consists of
the sums of the rows and of the columns is sufficient, a fact that has
already been established in Andersen (1977, 1980) or on p. 41 in Fischer and
Molenaar (1995) for the standard Rasch model and extended to the polytomous
Rasch model in Andrich (2010). But the above representation shows that one
can reduce the mixture Rasch model further to the statistic $%
(T_{1},\ldots,T_{m},N_{0},\ldots,N_{m})$ in a second step. As $%
\sum_{k=1}^{m}T_{k}=\sum_{k=0}^{m}N_{k}=n$ \ we may remove two
components of $(T_{1},\ldots,T_{m},N_{0},\ldots,N_{m})$ without losing sufficiency of the statistic. Especially the statistic $%
(T,N)=(T_{1},\ldots,T_{m-1},N_{1},\ldots,N_{m})$ is sufficient and takes its values
in $\{0,\ldots,n\}^{2m-1}.$ Denoting the distribution of $(T,N)$ under $%
P_{\theta ,F}^{A}$ by $P_{\theta ,F}^{B}$ we arrive at the model%
\begin{equation} \label{B}
\mathcal{B}_{n,m}:=\left( \{0,\ldots,n\}^{2m-1},\mathcal{P}(\{0,\ldots,n\}^{2m-1}),(P_{\theta ,F}^{B})_{(\theta,F)\in \Theta \times {\cal F}}\right) .  
\end{equation}%
As explained in Section \ref{aseq}, sufficiency implies equivalence in Le Cam's sense so that we obtain the following statement.

\begin{theo}
\label{T:0} The experiments $\mathcal{A}_{n,m}$ in (\ref{A}) and $\mathcal{B}%
_{n,m}$ in (\ref{B}) are equivalent.
\end{theo}

Put $\left\langle \theta ,b\right\rangle =\sum\nolimits_{i=1}^{m}\theta
_{i}b_{i},b \in \{0,1\}^{m}$ and
$$ \mathbb{S}(k,m) = \Big\{ b \in \{0,1\}^m : \sum\nolimits_{i=1}^{m} b_{i} = k \Big\}\,. $$
To study the distribution of $Y_j$ on $\{0,1\}^m = \bigcup_{k=0}^m \mathbb{S}(k,m)$ we deduce from (\ref{A}) that
\begin{equation}
P_{\theta ,F}^{A}(Y_{j}=b)=\exp \left\{ -\left\langle \theta ,b\right\rangle
+G(k,\theta ,F)\right\}\,, \qquad b\in \mathbb{S}(k,m)\,.
\end{equation}%
Moreover, the $Y_{j}$ are i.i.d. which implies that $(N_{0},\ldots,N_{m})$ has the multinomial distribution $%
\mathsf{M}_{n,m,\theta ,F}$ with the cell probabilities%
\begin{eqnarray} \nonumber 
q_{k}(\theta ,F) &=&\,\sum_{b\in \mathbb{S}(k,m)}\int \Big\{ \prod_{i=1}^{m}%
\frac{\exp \{b_{i}(\beta -\theta_{i})\}}{1+\exp \{\beta -\theta_{i}\}}%
\Big\} dF(\beta ) \\ \label{eq:q}
&=&\sum_{b\in \mathbb{S}(k,m)}\exp \left\{ -\left\langle \theta
,b\right\rangle +G(k,\theta ,F)\right\}.
\end{eqnarray}%
The conditional distribution $\Gamma _{\theta ,F}(\cdot |i)$ of $Y_{j}$
given $S_{j}=i$ has the probability mass function
\begin{align*} \Gamma
_{\theta ,F}(\{b\}|i) \, := \, 
P_{\theta ,F}^{A}(Y_{j}=b|S_{j}=i) & =\frac{\exp \left\{ -\left\langle \theta
,b\right\rangle \right\} }{\sum_{c\in \mathbb{S}(i,m)}\exp \left\{
-\left\langle \theta ,c\right\rangle \right\} } \cdot 1_{\mathbb{S}(i,m)}(b) \\
& = \frac{\exp\big(-\sum_{k=1}^{m-1} \vartheta_k b_k\big)}{\sum_{c \in \mathbb{S}(i,m)} \exp\big(-\sum_{k=1}^{m-1} \vartheta_k c_k\big)} \cdot 1_{\mathbb{S}(i,m)}(b)\,,
\end{align*}
where $\vartheta_k := \theta_k - \theta_m$. Writing $Y_j^{m-1} = (X_{j,1},\ldots,X_{j,m-1})$ the event $\{Y_j^{m-1} = b\}$ equals the union of $\{Y_j = (b,0)\}$ and $\{Y_j = (b,1)\}$ for any $b\in \{0,1\}^{m-1}$ so that
\begin{align} \label{eq:Utheta}
P_{\theta,F}^A\big(Y_j^{m-1} = b \mid S_j = i\big) & = \frac{\exp\big(-\sum_{k=1}^{m-1} \vartheta_k b_k\big)}{\sum_{c \in \mathbb{S}(i,m)} \exp\big(-\sum_{k=1}^{m-1} \vartheta_k c_k\big)} \cdot 1_{\mathbb{B}(i,m)}(b)\,, 
\end{align}
where $\mathbb{B}(i,m) := \mathbb{S}(i-1,m-1) \cup \mathbb{S}(i,m-1)$. The conditional measure of $Y_j^{m-1}$ given $S_j=i$ under $P_{\theta,F}^A$ is denoted by $\mathsf{U}_{\vartheta,F}(\cdot\mid i)$. As the random vectors $(Y_j^{m-1},S_j)$, $j=1,\ldots,n$, are independent the conditional measure of $T$ in (\ref{eq:neu}) given $S_1,\ldots,S_n$ under $P_{\theta,f}^A$ turns out to be
$$ {\cal L}(T\mid S_1,\ldots,S_n) \, = \, *_{j=1}^n \mathsf{U}_{\vartheta,F}(\cdot \mid S_j) \, = \, *_{k=1}^m \mathsf{U}_{\vartheta,F}^{*,N_k}(\cdot\mid k)\,, $$
where $*$ denotes convolution. Therein we have used that convolution is a commutative operation and that $\mathsf{U}_{\vartheta ,F}(\cdot |0)=\delta _{0}$. Since the random measure ${\cal L}(T\mid S_1,\ldots,S_n)$ is measurable in the $\sigma$-field generated by $N$ we conclude that
$$ {\cal L}(T\mid N) \, = \, *_{k=1}^m \mathsf{U}_{\vartheta,F}^{*,N_k}(\cdot \mid k)\,. $$
This proves
\begin{theo} \label{T:dis}
For the observation $(T,N) = (T_1,\ldots,T_{m-1},N_1,\ldots,N_m)$ in the experiment $\mathcal{B}_{n,m}$ in (\ref{B}), the random vector $\left(N_0:=
n-\sum_{i=1}^{m}N_{i},N_{1},\ldots,N_{m}\right) $ has a multinomial
distribution with the cell probabilities $q_{k}(\theta ,F)$; and $*_{i=1}^m \mathsf{U}_{\vartheta,F}^{*,N_i}(\cdot \mid i)$ is the conditional distribution of $T$ given $N$.
\end{theo}

It is remarkable that the conditional distribution of $T$ given $N$ does not depend on the ability distribution $F$ but only on the difficulty parameter $\theta$. This fact has also been mentioned e.g. in Pfanzagl (1993) and Strasser (2012a,b). 

\section{High-dimensional Gaussian Approximation} \label{Gauss}

In this section we establish a general result on the approximation of the sum of high-dimensional independent binary random vectors by Gaussian models. Later we will apply this finding to the experiment ${\cal B}_{n,m}$. The results of Carter (2002), which are restricted to multinomial experiments, are included in a special setting. In particular those results are not applicable to the statistic $T$ in the experiment ${\cal B}_{n,m}$. Moreover we use a completely different strategy of proofs. 

The starting point of this section is a triangular array of independent
binary vectors $Y_{i,n}=(X_{1,i,n},\ldots,X_{d,i,n})$ where the dimension $d=d_n$ is allowed to tend to infinity moderately with respect to $n$. That rate will be made precise later. We write $W_n := y_0 + \sum_{i=1}^n Y_{i,n}$ for any deterministic $y_0 \in \mathbb{Z}^d$. As $W_n$ is a discrete random vector which takes its values in $\mathbb{Z}^d$ one cannot approximate the measure $P_{W_n}$ of $W_n$ by a continuous probability measure such as a normal distribution in the total variation sense. Therefore one has to apply a smoothing procedure to $W_n$. Concretely, a $d$-dimensional random vector $U$ is generated independently of $W_n$ and, then, $W_n$ and $U$ are added so that we consider the continuous probability measure ${\cal L}(W_n + U) = {\cal L}(W_n) * {\cal L}(U)$. 

Now suppose that $W_n$ represents the observation in a statistical experiment. Then the Markov kernel $K(x,\cdot) := {\cal L}(U+x)$ transforms ${\cal L}(W_n)$ into ${\cal L}(W_n+U)$. As an attempt for the inverse transformation, one could round each component of $W_n+U$ and denote the outcome by $[W_n+U]$. Carter (2002) applies this strategy where $U$ is uniformly distributed on the cube $[-1/2,1/2]^d$. Then $[W_n+U] = W_n$ so that the original data are reconstructed by the rounding procedure. In this case the experiment in which one observes $W_n$ is equivalent to the experiment in which the observation is $W_n+U$. 

It turns out that, in the experiment ${\cal B}_{n,m}$, the approach which involves uniformly distributed $U$ would require $d_n$ to increase only at a logarithmic rate in $n$ in order to obtain asymptotic equivalence to a Gaussian model. Therefore we consider ${\cal L}(U) = \mathsf{N}(0,b_n I)$ where $I$ denotes the $d\times d$-identity matrix and the sequence $(b_n)_n$ is allowed to tend to infinity. Now the random vector $W_n$ cannot be identified from $W_n+U$ but we will show that the total variation distance between ${\cal L}(W_n)$ and ${\cal L}([W_n+U])$ still tends to zero (uniformly with respect to the parameter) under some constraints so that the experiment in which one observes $W_n$ is asymptotically equivalent to the experiment which describes the observation of $W_n+U$.
    
We introduce the notation
$$ Y_{i,n}^{-j} := (X_{1,i,n},\ldots,X_{j-1,i,n},X_{j+1,i,n},\ldots,X_{d,i,n})\,, $$ 
and
\begin{align*} p_{j,i} & := E(X_{j,i,n} \mid Y_{i,n}^{-j})\,, \qquad 
\mu_j  := \sum_{i=1}^n p_{j,i}\,, \qquad
\sigma^2_j  := \sum_{i=1}^n p_{j,i} (1 - p_{j,i})\,.
\end{align*} 
Moreover we define
\begin{equation} \label{eq:kappa} \kappa \, := \, \inf_{\theta'\in\Theta'} \min_{i=1,\ldots,n} \min_{j = 1,\ldots,d}\, E p_{j,i} (1-p_{j,i})\,, \end{equation}
when the distributions of the $Y_{i,n}$ are indexed by a parameter $\theta'\in\Theta'$. In order to show asymptotic proximity between ${\cal L}(W_n)$ and its shifted versions we provide the following lemma.
\begin{lem} \label{L:4.1}
Fix any $\delta>0$ such that $\kappa > n^{-1/2 + \delta}$ and $n\kappa > 2$. Then the total variation distance between ${\cal L}(W_n)$ and ${\cal L}(W_n+l)$, for some deterministic $l\in \mathbb{Z}^d$, obeys the following upper bound
$$ \mathsf{TV}({\cal L}(W_n),{\cal L}(W_n+l)) \, \leq \, A \{\log(n\kappa)\}^{1/2} n^{-1/2} \kappa^{-1/2} \, \sum_{j=1}^d |l_j|\,, $$
for a universal constant $A\in (0,\infty)$. 
\end{lem} 
Lemma \ref{L:4.1} represents a robustness property of ${\cal L}(W_n)$ with respect to shifting the measure on the $\mathbb{Z}^d$-grid. That provides the major tool for the announced upper bound on the total variation distance between ${\cal L}(W_n)$ and ${\cal L}([W_n+U])$. 
\begin{lem} \label{L:4.2} Under the conditions of Lemma \ref{L:4.1} we have, for ${\cal L}(U) = \mathsf{N}(0,b_{n}I)$, $b_n>0$, that
\begin{equation*}
\mathsf{TV}\big(\mathcal{L}(W_{n}),\mathcal{L}([W_{n}+U])\big) \, \leq \,  A \{\log(n\kappa)\}^{1/2} n^{-1/2} \kappa^{-1/2} \, d_n \, (1/2 + b_n^{1/2})\,.
\end{equation*}
with $A$ as in Lemma \ref{L:4.1}. 
\end{lem}
Thus, if the right hand side of the inequality in Lemma \ref{L:4.2} tends to zero (uniformly with respect to a family of admitted measures of the $Y_{i,n}$, $i=1,\ldots,n$), the observation of $W_n$, on the one hand, and of $W_n+U$, on the other hand, represent asymptotically equivalent experiments. 

In the next step we will approximate the smooth distribution ${\cal L}(W_n+U)$ by the normal distribution whose expectation vector and covariance matrix coincide with those of $W_n+U$. We establish a central limit theorem (CLT) for independent binary random vectors with increasing dimension in the total variation sense. We write $\overline{\mu}_i$ and $\overline{\Lambda}_i$ for the expectation vector and the covariance matrix of $Y_{i,n}$, respectively. Accordingly, $\overline{\mu} = y_0 + \sum_{i=1}^n \overline{\mu}_i$ and $\overline{\Lambda} \, = \, \sum_{i=1}^n \overline{\Lambda}_i$ are the corresponding quantities of $W_n$. Preparatory to this CLT we provide a positive lower bound on the eigenvalues of partial sums of the matrices $\overline{\Lambda}_i$. 
\begin{lem} \label{L:4.3}
All eigenvalues of the matrix $\sum_{i\in {\cal N}} \overline{\Lambda}_i$, for any ${\cal N} \subseteq \{1,\ldots,n\}$, are bounded from below by $(\# {\cal N})\, \kappa / d$. 
\end{lem}
Besides Lemma \ref{L:4.3} also yields invertibility of the matrix $\overline{\Lambda}$ whenever $\kappa>0$. Another important result which will be used to derive the CLT is the asymptotic proximity of the smoothed version of each ${\cal L}(Y_{i,n})$ (i.e. ${\cal L}(Y_{i,n})$ convolved with some normal distribution $\mathsf{N}(0,\tilde{\Lambda})$) and the normal distribution with the same expectation vector and covariance matrix as ${\cal L}(Y_{i,n})*\mathsf{N}(0,\tilde{\Lambda})$. We provide
\begin{lem} \label{L:4.4}
Let $\tilde{\Lambda}$ be some positive definite $d\times d$-matrix. Then,
$$ \mathsf{TV}\big({\cal L}(Y_{i,n})*\mathsf{N}(0,\tilde{\Lambda}),\mathsf{N}(\overline{\mu}_i,\overline{\Lambda}_i + \tilde{\Lambda})\big) \, \leq \, B\, \lambda^{-3/2}\, d_n^3\,, $$ 
for a universal constant $B\in (0,\infty)$ where $\lambda$ denotes the smallest eigenvalue of the matrix $\tilde{\Lambda}$.
\end{lem} 
We are now ready to prove a strong CLT for sums of independent binary random vectors.
\begin{lem} \label{L:4.5}
If $(b_n)_n$ is bounded away from zero then
$$ \mathsf{TV}\big({\cal L}(W_{n})*\mathsf{N}(0,b_n I),\mathsf{N}(\overline{\mu},\overline{\Lambda} + b_n I)\big) \, \leq \, C\, b_n^{-1/2}\, d_n^4 / \kappa_n\,, $$
for $\kappa = \kappa_n$ with a universal constant $C$. 
\end{lem}
Now we have a fully Gaussian random variable with the law $\mathsf{N}(\overline{\mu},\overline{\Lambda} + b_n I)$ where $\overline{\mu}$ and $\overline{\Lambda}$ represent the expectation vector and the covariance matrix of the original random vector $W_n$. Therefore the term $b_n I$ should be removed in the covariance matrix of the new random vector. By a famous formula which governs the Hellinger distance between normal distributions we deduce 
\begin{lem} \label{L:4.6}
We have that
$$ \mathsf{TV}\big(\mathsf{N}(\overline{\mu}+,\overline{\Lambda} + b_n I),\mathsf{N}(\overline{\mu},\overline{\Lambda})\big) \, \leq \, 2\sqrt{2}\, b_n\, n^{-1}\, \kappa_n^{-1}\, d_n^{3/2}\,. $$
\end{lem}

Piecing together the Lemmata \ref{L:4.2}, \ref{L:4.5} and \ref{L:4.6}, we derive the following central theorem which allows to approximate statistical experiments, in which one observes a sum of independent binary random vectors, by Gaussian experiments. Assume that the distributions of $Y_{i,n}$, $i=1,\ldots,n$, and $y_0\in \mathbb{Z}^d$ are indexed by a parameter $\theta'$, which lies in a set $\Theta'$. Then the experiment ${\cal X}_{n}$ describes the observation of the random vector $W_n$. Furthermore we define the Gaussian experiment ${\cal Z}_{n}$ by
$$ {\cal Z}_{n} \, := \, \big(\mathbb{R}^d,\mathfrak{B}(\mathbb{R}^d),\{\mathsf{N}(\overline{\mu}_{\theta'},\overline{\Lambda}_{\theta'})\}_{\theta'\in\Theta'}\big)\,. $$   
The above consideration leads to the following theorem, which is one of our main results. 
\begin{theo} \label{T:3}
Suppose that $\kappa = \kappa_n > n^{-1/2+\delta}$ for some fixed $\delta>0$ and $n$ sufficiently large; and that $\inf (b_n)_n > 0$. Then the Le Cam distance between the experiments ${\cal X}_n$ and ${\cal Z}_n$ satisfies
$$ \Delta\big({\cal X}_n,{\cal Z}_n\big) \, \leq \, \mbox{const.}\cdot \big(\{\log(n\kappa_n)\}^{1/2} n^{-1/2} \kappa_n^{-1/2} \, d_n \, b_n^{1/2} \, + \, b_n^{-1/2}\, d_n^4 / \kappa_n \, + \,  b_n\, n^{-1}\, \kappa_n^{-1}\, d_n^{3/2}\big)\,, $$
for some universal constant.   \end{theo}

\begin{rem} 
The Markov kernel which transforms ${\cal X}_n$ into ${\cal Z}_n$ in Theorem \ref{T:3} equals $x \mapsto \mathsf{N}(x,b_n I)$, $x \in \mathbb{Z}^d$; and the inverse transformation is carried out by rounding each component of the observation from ${\cal Z}_n$.
\end{rem} 

Pointing out the dominating terms, the upper bound on the Le Cam distance which is provided in Theorem \ref{T:3} converges to zero as $n\to\infty$ whenever
\begin{equation} \label{eq:T.3}
\lim_{n\to\infty} \{\log(n\kappa_n)\} n^{-1} \kappa_n^{-1} \, d_n^2 \, b_n \, + \, b_n^{-1}\, d_n^8 \kappa_n^{-2} \, = \, 0\,.
\end{equation}

\section{Gaussian Approximation of the Mixture Rasch Model} \label{5}

In this section we apply the general Gaussianization scheme provided in Section \ref{Gauss} and, in particular, in Theorem \ref{T:3} to the experiment ${\cal B}_{n,m}$ in (\ref{B}). Therein we distinguish between the statistics $T$ and $N$. Obviously $d_n$ from Theorem \ref{T:3} equals $m-1$ and $m$ for the statistic $T$ and $N$, respectively, while the quantity $\kappa_n$ has to be studied in both settings. 

\subsection{Gaussian Model for the Difficulty Parameters}

The new statistical experiment, which is denoted by ${\cal C}_{n,m}$, describes the observation of $(T^*,N)$ where $N$ is as in the experiment ${\cal B}_{n,m}$. Let $T^*$ be an $(m-1)$-dimensional random vector whose conditional distribution given $N$ is $\mathsf{N}\big(E_{\theta,F}^B(T|N),\mbox{cov}_{\theta,F}^B(T|N)\big)$. We define the experiments
\begin{align*}
{\cal B}_{n,m}^{n'} &\, := \, \big(\{0,\ldots,n\}^{m-1},{\cal P}(\{0,\ldots,n\}^{m-1}), ({\cal L}_{\theta,F}^B(T\mid N=n'))_{\theta,F}\big)\,, \\
{\cal C}_{n,m}^{n'}  &\, := \, \big(\mathbb{R}^{m-1},\mathfrak{B}(\mathbb{R}^{m-1}),({\cal L}_{\theta,F}^B(T^*\mid N=n'))_{\theta,F}\big)\,.
\end{align*}
Now we consider sequences of experiments indexed by the random vector $N$. Note that
\begin{equation} \label{eq:5.1} \Delta({\cal B}_{n,m},{\cal C}_{n,m}) \, \leq \, \sup_{\theta,F}\, E_{\theta,F}^B\, \Delta({\cal B}_{n,m}^{N},{\cal C}_{n,m}^{N})\,. \end{equation}
By Theorem \ref{T:dis} the observation in the experiment ${\cal B}_{n,m}^{n'}$ can be written as the sum of $n'_1+\cdots+n'_{m-1}$ independent binary random vectors so that it has the structure of the random vector $W_n$ from Section \ref{Gauss} when putting $y_0 = n'_m \cdot (1,\ldots,1)$. The following lemma gives us a lower bound on $\kappa$ in (\ref{eq:kappa}).
\begin{lem} \label{L:5.1}
Assuming that ${\cal L}(W_{n}-y_0 \mid N) = *_{k=1}^{m-1} \mathsf{U}^{*,N_k}_{\vartheta,F}(\cdot \mid k)$ in the notation of Section \ref{Gauss}; and that $m\geq 3$, the quantity $\kappa$ in (\ref{eq:kappa}) satisfies 
$$ \kappa \, \geq \, \frac{\exp(-6R)}{(m-1)(1+\exp(2R))}\,. $$
\end{lem}  
Note that the number of $Y_{i,n}$, which is denoted by $n$ in Section \ref{Gauss}, equals $n'_1+\cdots+n'_{m-1}$ in the experiment ${\cal B}_{n,m}^{n'}$. Therefore the following assumption and lemma are required. We impose that every distribution $F$ in ${\cal F}$ has a Lebesgue density $f$; and that there exists an envelopping function $\overline{f}$ with $\int \overline{f}(x) dx < \infty$ such that
\begin{equation} \label{eq:cond1}
f \, \leq \, \overline{f} \mbox{ a.e.}, \qquad \forall F \in {\cal F}\,.
\end{equation}
Condition (\ref{eq:cond1}) represents a tightness property of ${\cal F}$. Then
\begin{lem} \label{L:5.2}
Under the conditions (\ref{eq:cond1}), $m\geq 3$ and $\rho\in (0,1)$ sufficiently large, we have that
$$ \lim_{n\to\infty} \sup_{\theta\in\Theta,F\in {\cal F}} P_{\theta,F}^B\big(N_0 + N_m > \rho n\big) \, = \, 0\,. $$
\end{lem} 
By (\ref{eq:5.1}) we deduce for some $\rho\in (0,1)$ from Lemma \ref{L:5.2} that
\begin{align*} \Delta({\cal B}_{n,m},{\cal C}_{n,m}) & \, \leq \, 2  \sup_{\theta,F} P_{\theta,F}^B\big(N_0 + N_m > \rho n\big) \, + \, \sup_{\theta,F}\, E_{\theta,F}^B\, 1_{[(1-\rho)n,\infty)}\Big(\sum_{j=1}^{m-1} N_j\Big) \cdot \Delta({\cal B}_{n,m}^{N},{\cal C}_{n,m}^{N})\,, \end{align*}
where the latter term tends to zero as $n\to\infty$ whenever 
$$ \lim_{n\to\infty} n^{-1} m^3 b_n + b_n^{-1} m^{10}\, = \, 0\,, $$
thanks to Lemma \ref{L:5.1}, Theorem \ref{T:3} and equation (\ref{eq:T.3}). The convergence of the first term is guaranteed by Lemma \ref{L:5.2}. We establish asymptotic equivalence between the experiments ${\cal B}_{n,m}$ and ${\cal C}_{n,m}$ under some constraints. 
\begin{theo} \label{T:4}
Assume (\ref{eq:cond1}); $m=m_n\geq 3$; that there is some $\beta>13$ such that $\sup_n m_n^\beta / n < \infty$. Then the selection $b_n \asymp n^\alpha$ with $\alpha \in (10/\beta,1-3/\beta)$ yields asymptotic equivalence of the experiments ${\cal B}_{n,m}$ and ${\cal C}_{n,m}$ as $n\to\infty$. 
\end{theo}

Let us consider the conditional Gaussian distribution of the statistic $T^*$ given $N$ in the experiment ${\cal C}_{n,m}$. Since ${\cal L}(T\mid N) = *_{k=0}^m \mathsf{U}_{\vartheta,F}^{*,N_k}(\cdot\mid N)$ with $(T,N)$ as in the experiment ${\cal B}_{n,m}$ we have that $E_{\theta,F}^B(T\mid N) = - \nabla \Psi_N(\vartheta)$ and $\mbox{cov}_{\theta,F}^B(T\mid N) = \Delta \Psi_N(\vartheta)$ where
$$ \Psi_{n'}(\vartheta) \, := \, \sum_{k=0}^m n'_k \cdot \log\Big( \sum_{b\in \mathbb{S}(k,m)} \exp\Big\{- \sum_{l=1}^{m-1} \vartheta_l b_l\Big\}\Big), \qquad \vartheta = (\vartheta_1,\ldots,\vartheta_{m-1})\,, $$
and $\nabla$ and $\Delta$ denote the gradient and the Hessian matrix, respectively. 

We introduce the experiment ${\cal D}_{n,m}$ by
$$ {\cal D}_{n,m} \, := \, \big(\mathbb{R}^{2m-1},\mathfrak{B}(\mathbb{R}^{2m-1}),\big({\cal L}_{\theta,F}(N,T^{**})\big)_{\theta\in\Theta,F\in {\cal F}}\big)\,, $$
where $N$ is as in the experiment ${\cal C}_{n,m}$ and the conditional distribution of $T^{**}$ given $N$ equals $\mathsf{N}(\vartheta, \{\Delta \Psi_N(\vartheta)\}^{-1})$ if $N_0+N_m<n$; otherwise put $T^{**}=0$. By the Lemmata \ref{L:4.3} and \ref{L:5.1}, the matrix $\Delta \Psi_N(x)$ is invertible for all $x\in\mathbb{R}^{m-1}$ on the event $\{N_0+N_m < n\}$. Therein note that, for any $x\in\mathbb{R}^{m-1}$, there exist some $R>0$ and $\theta\in\Theta$ such that $\vartheta=x$. That also implies injectivity of the mapping $x\mapsto \nabla \Psi_N(x)$ on the domain $\mathbb{R}^{m-1}$ in the case of $N_0+N_m<n$. Now define the function $\Phi$ by
$$ \Phi(x,n') : \, = \, \begin{cases} (\nabla \Psi_{n'}(x),n') \,, & \mbox{ if }n'_1+\cdots+n'_{m-1} \neq 0\,, \\
(0,n')\,, & \mbox{ otherwise.} \end{cases} $$
By ${\cal D}'_{n,m}$ we define the experiment in which one observes $\Phi(T^{**},N)$ with $(T^{**},N)$ as in ${\cal D}_{n,m}$. Clearly $N$ is uniquely reconstructable from $\Phi(T^{**},N)$. If $N_1+\cdots+N_{m-1}=0$ then $T^{**}=0$; otherwise the injectivity of $x \mapsto \nabla \Psi_N(x)$ enables us to identify $T^{**}$. Therefore the experiments ${\cal D}_{n,m}$ and ${\cal D}'_{n,m}$ are equivalent in Le Cam's sense. Then its suffices to establish that
\begin{equation} \label{eq:Psi}
\lim_{n\to\infty} \sup_{\theta\in\Theta,F\in {\cal F}} \mathsf{TV}\big({\cal L}_{\theta,F}^D(\Phi),P_{\theta,F}^C\big) \, = \, 0\,, \end{equation}
in order to show the following theorem.
\begin{theo} \label{T:2}
Under the conditions of Theorem \ref{T:4} the experiments ${\cal C}_{n,m}$ and ${\cal D}_{n,m}$ are asymptotically equivalent as $n\to\infty$.
\end{theo}
The experiment ${\cal D}_{n,m}$ has the advantage compared to ${\cal C}_{n,m}$ that the directly observed statistic $T^{**}$ represents an asymptotically unbiased estimator of $\vartheta$. This will be exploited in Section \ref{Conf}.

\subsection{Gaussian Model for the Ability Distribution}

We focus on the multinomial statistic $N$ in the experiment ${\cal D}_{n,m}$. If we can show that the sub-experiment in which only $N$ is observed is asymptotically equivalent to the experiment which describes the observation of $N^*$ with
$$ {\cal L}(N^*) \, = \, \mathsf{N}\big(E_{\theta,F}^D N , \mbox{cov}_{\theta,F}^D(N)\big)\,, $$
then we have asymptotic equivalence of ${\cal D}_{n,m}$ and the experiment ${\cal E}_{n,m}$ which is defined by
$$ {\cal E}_{n,m} \, := \, \big(\mathbb{R}^{2m-1}, \mathfrak{B}(\mathbb{R}^{2m-1}) , ({\cal L}_{\theta,F}(T^{**},N^*))_{\theta\in \Theta,F\in {\cal F}}\big)\,, $$
such that ${\cal L}_{\theta,F}^E(T^{**}\mid N^*) = \mathsf{N}(\vartheta,\{\Delta \Psi_{[N^*]_+}(\vartheta)\}^{-1})$ if $[N_1^*]_+ + \cdots + [N_{m-1}^*]_+ > 0$ (put $T^{**}:=0$ otherwise) where $[x]_+$ denotes $(\max\{[x_j],0\})_{j=1,\ldots,m}$ for any $x\in \mathbb{R}^m$. Note that, for all $\theta\in\Theta$ and $F\in {\cal F}$, we have that 
$$ {\cal L}_{\theta,F}^E(T^{**}\mid N^*=N) = {\cal L}_{\theta,F}^D(T^{**}\mid N)\,, \quad \mbox{ a.s.}\,, $$
for $N$ as in the experiment ${\cal D}_{n,m}$. Moreover, by the multinomial distribution of $N$, we immediately derive that
\begin{align*}
E_{\theta,F}^D N & \, = \, n\cdot \tilde{q}(\theta,F) \, := \, n\cdot (q_1(\theta,F),\ldots,q_m(\theta,F))^T\,, \\
\mbox{cov}_{\theta,F}^D(N)  & \, = \, n \cdot \big(\tilde{Q}(\theta,F) - \tilde{q}(\theta,F) \tilde{q}(\theta,F)^T\big)\,, \end{align*}
where $q_k(\theta,F)$ is as in (\ref{eq:q}) and $\tilde{Q}(\theta,F)$ denotes the $(m-1)\times (m-1)$-diagonal matrix whose $(k,k)$th entry equals $q_k(\theta,F)$. The asymptotic equivalence of ${\cal D}_{n,m}$ and ${\cal E}_{n,m}$ is shown by a direct application of Theorem \ref{T:3} where the quantity $\kappa$ in (\ref{eq:kappa}) has to be bounded from below again. Therefore a constraint on the tail behaviour of the Lebesgue density $f$ of the ability distribution $F$ is required; concretely we assume that 
\begin{equation} \label{eq:tail}
f(x) \, \geq \, D_0 \exp\big\{-D_1 |x|\big\}\,, \qquad \forall x\in \mathbb{R}, F\in {\cal F}\,,
\end{equation}
for some universal positive constants $D_0$ and $D_1$. As an alternative for condition (\ref{eq:tail}) we may consider $m=m_n$ as bounded with respect to $n$. Then Gaussian models for $F$ are still included. In the notation of Section \ref{Gauss} it holds that
$$ \kappa \, \geq \, \inf_{\theta,F} \inf_{k=1,\ldots,m} q_0(\theta,F)\, q_k(\theta,F) / \big(q_0(\theta,F) + q_k(\theta,F)\big)\,. $$
Thus a lower bound on the $q_k(\theta,F)$ is needed.
\begin{lem} \label{L:5.3}
Under condition (\ref{eq:tail}) we obtain that
$$ \inf_{\theta,F} \inf_{k=0,\ldots,m} q_k(\theta,F) \, \geq \, \mbox{const.}\cdot m^{-3/2-D_1}\,, $$
for a universal constant factor.
\end{lem}
Hence $\kappa_n \asymp m^{-3/2-D_1}$ so that, by (\ref{eq:T.3}), the following statement is evident.
\begin{theo} \label{T:5}
Assume the constraints of Theorem \ref{T:4}; condition (\ref{eq:tail}); and the existence of some $\beta > 3 D_1 + 29/2$ such that 
$\sup_n m_n^\beta / n < \infty$. Then the selection $b_n \asymp n^\alpha$ with $\alpha \in \big((11+2D_1)/\beta,1 - (7/2+D_1)/\beta\big)$ yields asymptotic equivalence of the experiments ${\cal D}_{n,m}$ and ${\cal E}_{n,m}$ as $n\to\infty$. 
\end{theo}
Instead of condition (\ref{eq:tail}) one can assume that $m=m_n$ is bounded in $n$ and the claim of Theorem \ref{T:5} remains valid. 

Thanks to the multinomial distribution of the statistic $N$ in the experiment ${\cal B}_{n,m}$ a transformation of the experiment ${\cal E}_{n,m}$ (in particular, of the statistic $N^*$) is possible in order to obtain independent components. Similar arguments have been used in Carter (2002). We introduce the $(m+1)$-dimensional random vector $N^{**}$ with ${\cal L}(N^{**}) = \mathsf{N}(n q(\theta,F),n Q(\theta,F))$ where
\begin{align*}
q(\theta,F) & \, := \, \big(q_0(\theta,F),\ldots,q_m(\theta,F)\big)^T\,, \\
Q(\theta,F) & \, := \, \big\{1_{\{j\}}(k) q_k(\theta,F)\big\}_{j,k=0,\ldots,m}\,. \end{align*}
Then the conditional distribution of $T^{**}$ given $N^{**}$ equals $\mathsf{N}(\vartheta,\{\Delta \Psi_{[\tau(N^{**})]_+}(\vartheta)\}^{-1})$ on the event $\{[\tau(N^{**})]_{+,1} + \cdots + [\tau(N^{**})]_{+,m} > 0\}$ (again $T^{**}:=0$ otherwise), where the function $\tau$ from $\mathbb{R}^{m+1}$ to $\mathbb{R}^m$ is defined by 
$$ \tau \, : \, x = (x_0,\ldots,x_m) \, \mapsto \, (x_1,\ldots,x_m) \cdot n / \max\big\{\zeta,\sum_{j=0}^m x_j\big\}\,, $$
 for some deterministic $\zeta>0$ still to be chosen. We consider the experiment
$$ {\cal F}_{n,m} \, := \, \big(\mathbb{R}^{2m},\mathfrak{B}(\mathbb{R}^{2m}),({\cal L}_{\theta,F}(T^{**},N^{**}))_{\theta,F}\big)\,. $$
In order to show asymptotic equivalence of ${\cal E}_{n,m}$ and ${\cal F}_{n,m}$ we consider the statistic $N^{**}$ from the experiment ${\cal F}_{n,m}$ and the sum of its components, which we call $V$. As $N^{**}$ can be uniquely reconstructed from $(\tau(N^{**}),V)$ we derive equivalence of ${\cal F}_{n,m}$ and the experiment ${\cal F}'_{n,m}$ in which $(T^{**},\tau(N^{**}),V)$ is observed. It holds that 
\begin{align*}
{\cal L}_{\theta,F}(N^{**} \mid V) & \, = \, \mathsf{N}\big(V q(\theta,F),n Q(\theta,F) - n q(\theta,F) q(\theta,F)^T\big)\,, \\
 {\cal L}_{\theta,F}(\tau(N^{**})\mid V) & \, = \, \mathsf{N}\big(n \tilde{q}(\theta,F) V / \max\{T,\zeta\}, n^3 \big(\tilde{Q}(\theta,F) - \tilde{q}(\theta,F) \tilde{q}(\theta,F)^T\big) / (\max\{V,\zeta\})^2\big)\,. \end{align*}
The following asymptotic approximation is required.
\begin{lem} \label{L:5.4}
Assume the conditions of Theorem \ref{T:5} and select $\zeta = n/2$. Then,
$$ \lim_{n\to\infty} \sup_{\theta,F} E_{\theta,F} \mathsf{TV}\big({\cal L}_{\theta,F}(\tau(N^{**})\mid V), \mathsf{N}(n \tilde{q}(\theta,F), n \tilde{Q}(\theta,F) - n \tilde{q}(\theta,F) \tilde{q}(\theta,F)^T)\big) \, = \, 0\,. $$
\end{lem}

As the conditional distribution of $T^{**}$ given $\tau(N^{**})$ and $V$ equals that given $\tau(N^{**})$, Lemma \ref{L:5.4} provides asymptotic equivalence of ${\cal F}'_{n,m}$ and ${\cal F}''_{n,m}$ where the latter experiment describes the observation of $(T^{**},N^{*},V^*)$ where $V^*$ and $V$ are identically distributed but $V^*$ and $(T^{**},N^{*})$ are independent. As, in addition, ${\cal L}_{\theta,F}(V^*) = \mathsf{N}(n,n)$, the distribution of $V^*$ does not depend on $\theta$ or $F$ and, thus, $V^*$ can be omitted without losing any information on $(\theta,F)$. Therefore, ${\cal F}''_{n,m}$ and ${\cal E}_{n,m}$ are equivalent so that the following result has been established.
\begin{theo} \label{T:6}
Under the conditions of Theorem \ref{T:5}, the experiments ${\cal E}_{n,m}$ and ${\cal F}_{n,m}$ are asymptotically equivalent as $n\to\infty$.
\end{theo}
 
\section{Applications} \label{Conf}

In this section we apply the Gaussian models of Section \ref{5}, which have now been proved to be asymptotically equivalent to the mixture Rasch model ${\cal A}_{n,m}$, to develop asymptotic inference. In particular we will construct an asymptotic confidence ellipsoid for the difficulty parameters. Thus the results carry over to the original mixture Rasch model. 

Let $T^{**}$ be the part of the observation from the experiment ${\cal D}_{n,m}$ where ${\cal L}(T^{**}\mid N) = \mathsf{N}\big(\vartheta,\{\Delta \Psi_{N}(\vartheta)\}^{-1}\big)$ if $N_1+\cdots+N_{m-1} > 0$. We define the random ellipsoid
$$ \hat{E} \, := \, \Big\{x\in \mathbb{R}^{m} \, : \, \sum_{j=1}^{m} x_j = 0\,, \, (x - Z T^{**})^T Z (Z^T Z)^{-1} \Delta \Psi_{N}(T^{**})  (Z^\dagger Z)^{-1} Z^\dagger (x - Z T^{**}) \, \leq \, \iota\Big\}\,, $$
for some threshold $\iota>0$ to be determined and the $m\times(m-1)$-matrix
$$ Z \, := \, \begin{pmatrix} 1 - 1/m & -1/m & \cdots & -1/m \\
                                -1/m & 1 - 1/m &\cdots & -1/m \\
                                 \vdots &\vdots & \vdots & \vdots \\
                                 -1/m & -1/m & \cdots & 1-1/m \\
                                 -1/m & -1/m & \cdots & -1/m \end{pmatrix}\,. $$
Note that $\theta = Z \vartheta$ for all $\theta\in\Theta$ thanks to the definition (\ref{eq:Theta}). In order to motivate the selection of $\iota$ we give an oracle version of $\hat{E}$ by
$$ \tilde{E} \, := \, \Big\{x\in \mathbb{R}^{m} \, : \, \sum_{j=1}^{m} x_j = 0\,, \, (x - Z T^{**})^\dagger Z (Z^\dagger Z)^{-1} \Delta \Psi_{N}(\vartheta)  (Z^\dagger Z)^{-1} Z^\dagger (x - Z T^{**}) \, \leq \, \iota\Big\}\,. $$
Concretely we have replaced $T^{**}$ by its conditional expectation given $N$ in the argument of $\Delta \Psi_{N}$. Conditionally on $N$ under $\{N_1+\cdots+N_{m-1}>0\}$ we may represent $T^{**}$ by $T^{**} = \vartheta + \{\Delta \Psi_{N}(\vartheta)\}^{-1/2} \varepsilon$ where here $\varepsilon$ denotes an $(m-1)$-dimensional random vector with independent $\mathsf{N}(0,1)$-distributed components. On this event it follows that
$$ P^D_{\theta,F}\big(\theta \in \tilde{E} \mid N\big) \, = \, P\big(\big|\varepsilon\big|^2 \leq \iota\big)\,. $$
That inspires us to choose $\iota$ as the $\alpha$-quantile of the $\chi^2(m-1)$-distribution, i.e. $\iota = F_{m-1}^{-1}(\alpha)$ where $F_{m-1}$ denotes the $\chi^2(m-1)$-distribution function, for some given $\alpha \in (0,1)$. Then, 
$$ \liminf_{n\to\infty} \inf_{\theta,F} P_{\theta,F}^D\big(\theta \in \tilde{E}\big) \, \geq \, \alpha\,, $$
as $P_{\theta,F}^D(N_0+N_m < n)$ tends to zero uniformly in $\theta$ and $F$. Focusing on the ellipsoid $\hat{E}$ we provide the following result.
\begin{theo} \label{T:appl1}
In the experiment ${\cal D}_{n,m}$ we have 
$$ \limsup_{n\to\infty} \sup_{\theta,F} P^D_{\theta,F}\big(\theta \not\in \hat{E}\big) \, \leq \, 1 - \alpha\,, $$
under the assumptions of Theorem \ref{T:4} for any fixed $\alpha\in (0,1)$. The maximal axis $\hat{e}$ of $\hat{E}$ obeys the following asymptotic upper bound
$$ \lim_{c\to\infty} \limsup_{n\to\infty} \sup_{\theta,F} P^D_{\theta,F}\big(\hat{e} > c \cdot m_n/\sqrt{n}\big) \, = \, 0\,. $$
\end{theo}

\begin{rem} 
Theorem \ref{T:appl1} shows that $\hat{E}$ is an asymptotic $\alpha$-confidence ellipsoid for $\theta$ in the experiment ${\cal D}_{n,m}$. The maximal axis of this ellipsoid shrinks to zero at the rate ${\cal O}_P(m_n/\sqrt{n})$ as $n\to\infty$. By the Theorems \ref{T:0}, \ref{T:dis}, \ref{T:4} and \ref{T:2} the properties established in Theorem \ref{T:appl1} extend to the original mixture Rasch model (experiment ${\cal A}_{n,m}$) after applying the appropriate Markov kernel which transforms experiment ${\cal D}_{n,m}$ to ${\cal A}_{n,m}$. Note that the asymptotic confidence region is uniform with respect to the parameter $\theta\in\Theta$ and $F\in {\cal F}$. Thus we have developed a stronger version of asymptotic confidence regions than in the usual setting where $\theta$ and $F$ are viewed as fixed, i.e. $\theta$ and $F$ must not change in $n$. This is thanks to the fact that we have used asymptotic approximation with respect to the Le Cam distance rather than central limit laws for specific estimators in terms of weak convergence, where the latter results are commomly used to construct asymptotic confidence regions.   
\end{rem}

\section{Conclusions and Outlook}

In this paper we derive asymptotically equivalent Gaussian experiments for the mixture Rasch model. In Section \ref{Conf}, asymptotic statistical inference on the difficulty parameters is provided based on these Gaussian experiments. But the asymptotic equivalence of the experiment ${\cal F}_{n,m}$ and the original mixture Rasch model ${\cal A}_{n,m}$ also opens the perspective for nonparametric inference on the ability distribution. While this goal exceeds the framework of the current paper the authors are working on this issue and intend to present their future results in a separate paper.

\section{Proofs}

\noindent {\it Proof of Lemma \ref{L:4.1}}: Thanks to the shift-invariance of the total variation distance we may put $y_0=0$ without any loss of generality. Note that we may write
$$ {\cal L}(W_n + l) \, = \, {\cal L}(W_n) * \big(*_{j=1}^d \delta_{l_j e_j}\big)\,, $$
where $e_j$ denotes the vector with its $j$th component equal to $1$ while all other components vanish. By a telescoping sum we deduce that
$$ \mathsf{TV}({\cal L}(W_n + l), {\cal L}(W_n)) \, \leq \, \sum_{j=1}^d \mathsf{TV}({\cal L}(W_n),{\cal L}(W_n) * \delta_{l_j e_j}) \, \leq \, \sum_{j=1}^d |l_j| \, \mathsf{TV}({\cal L}(W_n),{\cal L}(W_n) * \delta_{e_j})\,. $$
We have that
\begin{align} \nonumber \mathsf{TV}({\cal L}(W_n),{\cal L}(W_n) * \delta_{e_j}) & \, = \, \sum_{w \in \mathbb{Z}^d} \big|P(W_n=w) - P(W_n=w-\delta_{e_j})\big| \\ \label{eq:Pr.4.1.0} & \, \leq \, E \sum_{u \in \mathbb{Z}} \big|P(W_{n,j}=u\mid {\cal Y}_{n,j}) - P(W_{n,j}=u-1\mid {\cal Y}_{n,j})\big|\,, \end{align}
where ${\cal Y}_{n,j}$ denotes the $\sigma$-field generated by $Y_{1,n}^{-j},\ldots,Y_{n,n}^{-j}$. By Fourier inversion we obtain that
$$ P(W_{n,j} = u \mid {\cal Y}_{n,j}) \, = \, \frac1{2\pi} \int_{-\pi}^\pi \exp\{-iux\} \psi_{W_{n,j}\mid {\cal Y}_{n,j}}(x) dx\,, $$
for all $u\in \mathbb{Z}$ where $\psi_{Z}$ denotes the characteristic function of a random variable $Z$. Since
\begin{align*} |\psi_{W_{n,j}\mid {\cal Y}_{n,j}}(x)| & \, = \, \prod_{i=1}^n \big|\psi_{X_{j,i,n}\mid Y_{i,n}^{-j}}(x)\big| \, = \, \prod_{i=1}^n \big| \exp\{ix\} p_{j,i} + 1 - p_{j,i}\big| \\
& \, \leq \, \prod_{i=1}^n \exp\big\{ - 2 p_{j,i} (1-p_{j,i}) x^2 / \pi^2\big\} \, = \, \exp\big\{ - 2 \sigma_j^2 x^2 / \pi^2\big\}\,,
\end{align*}
for all $x\in [-\pi,\pi]$ it follows that
\begin{align*}
\big|P&(W_{n,j} = u \mid {\cal Y}_{n,j}) - P(W_{n,j} = u-1 \mid {\cal Y}_{n,j})\big| \\ & \, = \,  \frac1{2\pi} \Big| \int_{-\pi}^\pi \exp\{-iux\} \big(1 - \exp\{-ix\}\big) \psi_{W_{n,j}\mid {\cal Y}_{n,j}}(x) dx \Big| \\ &  \, \leq \, \frac1{\pi} \int_{0}^\pi x \exp\big\{ - 2 \sigma_j^2 x^2 / \pi^2\big\} dx \, \leq \, \pi / (4 \sigma_j^2)\,.
\end{align*}
Therefore the total variation distance between ${\cal L}(W_{n,j}\mid {\cal Y}_{n,j})$ and ${\cal L}(W_{n,j}+1\mid {\cal Y}_{n,j})$ is bounded from above by
\begin{align*}
& \sum_{u \in \mathbb{Z}} \big|P(W_{n,j}=u\mid {\cal Y}_{n,j}) - P(W_{n,j}=u-1\mid {\cal Y}_{n,j})\big| \\
& \, \leq \, \sum_{|u - \mu_j| \leq \tau \sigma_j + 1} \big|P(W_{n,j}=u\mid {\cal Y}_{n,j}) - P(W_{n,j}=u-1\mid {\cal Y}_{n,j})\big| \, + \, 2 \, P(|W_{n,j} - \mu_j| > \tau \sigma_j \mid {\cal Y}_{n,j}) \\
& \, \leq \, (2\tau \sigma_j + 3)  \pi^2 / (2 \sigma_j^2) + 4 \exp\{-\tau^2/4\} + 4 \exp\{-3 \sigma_j \tau / 4\}\,,
\end{align*} 
for any $\tau>0$ where Bernstein's inequality has been used in the last step. We introduce the event $A_j := \{\sigma^2_j > E \sigma_j^2 / 2\}$ where $E \sigma_j^2 \geq n \kappa$ and we put $\tau := \sqrt{c \cdot \log(n \kappa)}$ with a constant $c>0$ sufficiently large so that
\begin{align*}
E &\sum_{u \in \mathbb{Z}} \big|P(W_{n,j}=u\mid {\cal Y}_{n,j}) - P(W_{n,j}=u-1\mid {\cal Y}_{n,j})\big|   \\ & \hspace{7.5cm} \leq  2(1-P(A_j)) + A^* \{\log(n \kappa)\}^{1/2} n^{-1/2} \kappa^{-1/2}\,, \end{align*}
for a universal constant $A^* \in (0,\infty)$. Finally Hoeffding's inequality yields that
$$ 1 - P(A_j) \, \leq \, 2 \exp\big\{ - (E \sigma_j^2)^2 / (2 n)\big\} \, \leq \, 2 \exp\{ - n \kappa^2 / 2\}\,, $$
which completes the proof of the lemma. \hfill $\square$ \\

\noindent {\it Proof of Lemma \ref{L:4.2}}: As $W_n$ is $\mathbb{Z}^d$-valued it holds that $[W_n+U] = W_n + [U]$ so that
\begin{align*}
\mathsf{TV}({\cal L}([W_n+U]),{\cal L}(W_n)) & \, \leq \, \sum_{l\in \mathbb{Z}^d} \mathsf{TV}({\cal L}(W_n+l),{\cal L}(W_n)) \cdot P([U]=l) \\
& \, \leq \, A \{\log(n\kappa)\}^{1/2} n^{-1/2} \kappa^{-1/2} \, \sum_{j=1}^d E |[U_j]|\,,
\end{align*}
where $E|[U_j]| \, \leq \, 1/2 + E|U_j| \, \leq \, 1/2 + b_n^{1/2}$. \hfill $\square$ \\

\noindent {\it Proof of Lemma \ref{L:4.3}}: Let $\lambda$ be an arbitrary eigenvalue of the matrix $\sum_{i \in {\cal N}} \overline{\Lambda}_i$ with the corresponding unit eigenvector $v$. As $\overline{\Lambda}_i$ is the covariance matrix of $Y_{i,n}$ we deduce that
\begin{align*}
\lambda & \, = \, v^T \sum_{i \in {\cal N}} \overline{\Lambda}_i v \, = \, \sum_{i \in {\cal N}} v^T \overline{\Lambda}_i v \, = \, \sum_{i \in {\cal N}} \mbox{var}\Big(\sum_{k=1}^d v_k X_{k,i,n}\Big) \\ & \, = \, \sum_{i \in {\cal N}} E\, E\Big\{ \Big(\sum_{k=1}^d v_k (X_{k,i,n} - E X_{k,i,n})\Big)^2 \mid Y_{i,n}^{-l}\Big\} \, \geq \, \sum_{i \in {\cal N}} E\, \mbox{var}\big(v_l X_{k,l,n} \mid Y_{i,n}^{-l}\big) \, \geq \, v_l^2 \cdot (\# {\cal N}) \cdot \kappa\,,
\end{align*}
for all $l=1,\ldots,d$. Summing up both sides of the above inequality over $l=1,\ldots,d$ we obtain that
$$ d \lambda \, \geq \, (\# {\cal N}) \cdot \kappa\,, $$ 
which completes the proof. \hfill $\square$ \\

\noindent {\it Proof of Lemma \ref{L:4.4}}: Note that $\tilde{\Lambda} - \lambda I$ is a positive semi-definite matrix so that
$$ \mathsf{N}(0,\tilde{\Lambda}) \, = \, \mathsf{N}(0,\lambda I) * \mathsf{N}(0,\tilde{\Lambda} - \lambda I)\,, $$
from what follows that
$$ \mathsf{TV}\big({\cal L}(Y_{i,n})*\mathsf{N}(0,\tilde{\Lambda}),\mathsf{N}(\tilde{\mu}_i,\overline{\Lambda}_i + \tilde{\Lambda})\big) \, \leq \, \mathsf{TV}\big({\cal L}(Y_{i,n})*\mathsf{N}(0,\lambda I),\mathsf{N}(\tilde{\mu}_i,\overline{\Lambda}_i + \lambda I)\big)\,. $$
The distribution ${\cal L}(Y_{i,n})*\mathsf{N}(0,\lambda I)$ has the $d$-dimensional Lebesgue density
\begin{align*}
g_0(x) & \, = \, E \, (2\pi \lambda)^{-d/2} \exp\big\{ - |x - Y_{i,n}|^2 / (2\lambda)\big\}\,. \end{align*} 
Since $\mathsf{N}(\overline{\mu}_i,\overline{\Lambda}_i + \lambda I) = \mathsf{N}(\overline{\mu}_i,\overline{\Lambda}_i) * \mathsf{N}(0,\lambda I)$ the Lebesgue density $g_1$ of the distribution $\mathsf{N}(\overline{\mu}_i,\overline{\Lambda}_i + \lambda I)$ may be written as
$$ g_1(x) \, = \, E \, (2\pi\lambda)^{-d/2} \exp\big\{ - |x - Z_{i,n}|^2 / (2\lambda)\big\}\,, $$
where ${\cal L}(Z_{i,n}) = \mathsf{N}(\overline{\mu}_i,\overline{\Lambda}_i)$. The total variation distance between $\mathsf{N}(\overline{\mu}_i,\overline{\Lambda}_i + \lambda I)$ and ${\cal L}(Y_{i,n}) * \mathsf{N}(0,\lambda I)$ equals the $L_1(\mathbb{R}^d)$-distance between the densities $g_0$ and $g_1$. Thus,
\begin{align*}
\mathsf{TV}&\big({\cal L}(Y_{i,n})*\mathsf{N}(0,\tilde{\Lambda}),\mathsf{N}(\overline{\mu}_i,\overline{\Lambda}_i + \lambda I)\big) \\
& \, = \, (2\pi)^{-d/2} \, \int \big|E \exp\big\{ - |x - \lambda^{-1/2} Y_{i,n}|^2 / 2\big\} - E \exp\big\{ - |x - \lambda^{-1/2} Z_{i,n}|^2 / 2\big\}\big| dx\,.
\end{align*}
Taylor expansion around $x$ yields that
$$ \exp\big(-|x-\Delta|^2/2\big) \, = \, P_{2,x}(\Delta) + R_{2,x}(\Delta)\,, $$
for all $\Delta\in \mathbb{R}^d$ and any fixed $x\in \mathbb{R}^d$ where $P_{2,x}$ is a $d$-variate quadratic polynomial and $R_{2,x}$ is the corresponding remainder term. As the expectation vectors and the covariance matrices of $Y_{i,n}$ and $Z_{i,n}$ coincide we deduce that
$$ E \, P_{2,x}\big(\lambda^{-1/2} Y_{i,n}\big) \, = \, E\, P_{2,x}\big(\lambda^{-1/2} Z_{i,n}\big)\,. $$
Therefore,
\begin{align*}
\mathsf{TV}&\big({\cal L}(Y_{i,n})*\mathsf{N}(0,\lambda I),\mathsf{N}(\overline{\mu}_i,\overline{\Lambda}_i + \lambda I)\big) \\ &  \leq  (2\pi)^{-d/2} \Big( \int E \big|R_{2,x}\big(\lambda^{-1/2} Y_{i,n}\big)\big| dx + \int E \big|R_{2,x}\big(\lambda^{-1/2} Z_{i,n}\big)\big| dx \Big).
\end{align*}
Calculating the third order partial derivatives of $x \mapsto \exp(-|x|^2/2)$ we deduce that 
\begin{align} \nonumber (2\pi)^{-d/2} E \int \big|R_{2,x}(\lambda^{-1/2} Y_{i,n})\big| dx & \, \leq \, \lambda^{-3/2}\, E \Big(\sum_{k=1}^d |X_{k,i,n}|\Big)^3 \\
& \, \leq \, B^* \cdot \lambda^{-3/2} \cdot d^3\,, \end{align}
for some universal constant $B^* \in (0,\infty)$ where ${\cal L}(Z) = \mathsf{N}(0,1)$. When replacing $X_{k,i,n}$ by the $k$th component of $Z_{i,n}$ the identical upper bound applies (with a different constant $B^*$). Note that this component is $\mathsf{N}\big(E X_{k,i,n}, \mbox{var}\, X_{k,i,n}\big)$-distributed. \hfill $\square$ \\

\noindent {\it Proof of Lemma \ref{L:4.5}}: Again the shift-invariance of the total variation distance allows us to restrict to the case of $y_0=0$. As a telescoping sum, we consider 
\begin{align*}
\mathsf{TV}&\big({\cal L}(W_n) * \mathsf{N}(0,b_n I) , \mathsf{N}(\overline{\mu},\overline{\Lambda} + b_n I)\big) \\ & \, = \, \mathsf{TV}\big( \mathsf{N}(0,b_n I) * \big\{ *_{i=1}^n {\cal L}(Y_{i,n})\big\} , \mathsf{N}(0,b_n I) * \big\{ *_{i=1}^n \mathsf{N}(\overline{\mu}_i,\overline{\Lambda}_i)\big\}\big) \\
& \, \leq \, \sum_{k=1}^{n} \mathsf{TV}\big(  \big\{*_{i=k+1}^{n} {\cal L}(Y_{i,n})\big\} * \mathsf{N}(0,b_n I) * \big\{*_{i=1}^k \mathsf{N}(\overline{\mu}_i , \overline{\Lambda}_i)\big\}, \\ & \hspace{7cm} \big\{*_{i=k}^{n} {\cal L}(Y_{i,n})\big\} * \mathsf{N}(0,b_n I) * \big\{*_{i=1}^{k-1} \mathsf{N}(\overline{\mu}_i , \overline{\Lambda}_i)\big\}\big) \\
& \, \leq \, \sum_{k=1}^{n} \mathsf{TV}\big( \mathsf{N}( \overline{\mu}_k , \overline{\Lambda}_k + \tilde{\Lambda}_{k-1})  ,  {\cal L}(Y_{k,n}) *  \mathsf{N}( 0, \tilde{\Lambda}_{k-1})\big)\,, 
\end{align*}
where $\tilde{\Lambda}_{l} := b_n I + \sum_{i=1}^l \overline{\Lambda}_i$. By Lemma \ref{L:4.3} the smallest eigenvalue of $\tilde{\Lambda}_l$ is bounded from below by $b_n + l \kappa / d_n$. Then Lemma \ref{L:4.4} provides that
\begin{align*}
\mathsf{TV}& \big({\cal L}(W_n) * \mathsf{N}(0,b_n I) , \mathsf{N}(\overline{\mu},\overline{\Lambda} + b_n I)\big) \, \leq \, B \, d_n^3 \, \sum_{k=1}^n (b_n + (k-1) \kappa / d)^{-3/2} \\
& \, \leq \, B \, d_n^3 \, \Big( b_n^{-3/2} + \int_0^\infty (b_n + x \kappa / d_n)^{-3/2}dx \Big) \, = \, B \, d_n^3 b_n^{-3/2} + 2B\, d_n^4 b_n^{-1/2} / \kappa\,. \end{align*}
Thus the lemma has been shown. \hfill $\square$ \\ 
 
\noindent {\it Proof of Lemma \ref{L:4.6}}: The total variation distance between two distributions is bounded from above by twice their Hellinger distance. It follows from e.g. eq. (A.4) in Rei{\ss} (2011) that
\begin{align*}
\mathsf{TV}^2\big(\mathsf{N}(\overline{\mu},\overline{\Lambda} + b_n I),\mathsf{N}(\overline{\mu},\overline{\Lambda})\big) & \, \leq \, 8 b_n^2 \big\|\overline{\Lambda}^{-1}\big\|_F^2 \, = \, 8 b_n^2 \sum_{j=1}^{d_n} \lambda_j^{-2}\,, \end{align*}
where $\|\cdot\|_F$ denotes the Frobenius norm and $\lambda_j$, $j=1,\ldots,d_n$, are the eigenvalues of the matrix $\overline{\Lambda}$. Applying the lower bound on the eigenvalues provided in Lemma \ref{L:4.3} completes the proof of this lemma. \hfill $\square$ \\

\noindent {\it Proof of Lemma \ref{L:5.1}}: In the notation of Section \ref{Gauss} we assume some random vector $Y_{i,n} = (X_{1,i,n},\ldots,X_{d,i,n})$ with $d=m-1$ and ${\cal L}(Y_{i,n}) = \mathsf{U}_{\vartheta,F}(\cdot\mid k)$. Then
$$ p_{l,i} = E (X_{l,i,n} \mid Y_{i,n}^{-l}) = \frac{\exp(-\vartheta_l)}{1 + \exp(-\vartheta_l)} \cdot 1_{\{k-1\}}\Big(\sum_{q\neq l}^{m-1} X_{q,i,n}\Big) + 1_{\{k-2\}}\Big(\sum_{q\neq l}^{m-1} X_{q,i,n}\Big)\,, $$
so that
\begin{align*} E\, p_{l,i} (1-p_{l,i}) & \, = \, \frac{\exp(-\vartheta_l)}{(1 + \exp(-\vartheta_l))^2} \cdot \mathsf{U}_{\vartheta,F}\Big(\sum_{q\neq l}^{m-1} X_{q,i,n} = k-1\mid k\Big) \\
& \, = \, \frac{\exp(-\vartheta_l)}{1 + \exp(-\vartheta_l)} \cdot \sum_{b \in \mathbb{B}'(l,k,m)} \exp\Big(- \sum_{i\neq l}^{m-1} \vartheta_i b_i\Big) \, \big/ \, \sum_{c \in \mathbb{B}(k,m)} \exp\Big(- \sum_{i=1}^{m-1} \vartheta_i c_i\Big)\,,
\end{align*}
where $\mathbb{B}'(l,k,m)$ collects all $b \in \mathbb{B}(k,m)$ such that $\sum_{i\neq l}^{m-1} b_i = k-1$. Note that 
$$ \mathbb{B}(k,m) \, = \, \bigcup_{l=1}^{m-1} \mathbb{B}'(l,k,m)\,, $$
as $k \in \{1,\ldots,m-1\}$. Also we have
$$  \sum_{b \in \mathbb{B}'(l,k,m)} \exp\Big(- \sum_{i\neq l}^{m-1} \vartheta_i b_i\Big) \, \leq \, \exp(4R) \, \sum_{b \in \mathbb{B}'(l',k,m)} \exp\Big(- \sum_{i\neq l'}^{m-1} \vartheta_i b_i\Big)\,, $$
for all $l,l'\in \{1,\ldots,m-1\}$. It follows that
$$ E\, p_{l,i} (1-p_{l,i}) \, \geq \, \frac{\exp(-6R)}{(m-1)(1+\exp(2R))}\,, $$
which completes the proof. \hfill $\square$ \\

\noindent {\it Proof of Lemma \ref{L:5.2}}: Note that $N_0+N_m$ has a binomial distribution with the parameters $n$ and $q_0(\theta,F)+q_m(\theta,F)$ where $q_k(\theta,F)$ is defined in (\ref{eq:q}).  For any $s > 0$ we have 
\begin{align} \nonumber
q_0(\theta,F) & \, \leq \, \int_{x<-s} \overline{f}(x) dx \, + \, \big(1 + \exp(-s-R)\big)^{-m}\,, \\ \label{eq:Pr.Gauss.2.2}
q_1(\theta,F) & \, \leq \, \int_{x>s} \overline{f}(x) dx \, + \, \big(1 + \exp(-s-R)\big)^{-m}\,,
\end{align} 
for all $\theta\in\Theta$. For any fixed $\varepsilon>0$ we choose $s$ sufficiently large such that the first addends in both lines of (\ref{eq:Pr.Gauss.2.2}) are bounded from above by $\varepsilon/2$ and, then, $M$ sufficiently large such that for all $m>M$ the second addends in (\ref{eq:Pr.Gauss.2.2}) are smaller than $\varepsilon/2$. Thus, for all $m>M$, we obtain that $q_0(\theta,F)+q_m(\theta,F) < \varepsilon$. On the other hand, if $m\leq M$, we fix $s'>0$ sufficiently large such that
$$ \int_{|x|\geq s'} \overline{f}(x) dx \, < \, 1/2\,, $$
so that $\int_{|x|\leq s'} f(x) dx \, \geq \, 1/2$ holds true for all $F\in {\cal F}$. Then we consider the continuous positive mapping $T_m$, $m=2,\ldots,M$, with
$$ T_m(x,\theta) \, := \, \sum_{b \in \mathbb{S}(1,m)} \prod_{l=1}^m \frac{\exp(b_l [x-\theta_{l}])}{1 + \exp(x-\theta_{l})}\,, $$
which take its positive minimum on the compact domain $[-s',s']\times [-R,R]^m$. Hence $$ \inf_{\theta\in\Theta,F\in {\cal F}} \, q_1(\theta,F) \, > \, 0, \qquad \forall m=2,\ldots,M\,, $$ so that $\sup_{\theta\in\Theta,F\in {\cal F}} \, q_0(\theta,F)+q_m(\theta,F) < 1$. Thus we have shown that
$$ q := \sup_{\theta\in\Theta,F\in {\cal F}} \, \sup_{m\geq 3} \, q_0(\theta,F) + q_m(\theta,F) \, < \, 1\,. $$
Now we choose $\rho := (1+q)/2 \in (0,1)$ so that simple application of Chebyshev's inequality completes the proof. \hfill $\square$ \\
 
\noindent {\it Proof of Theorem \ref{T:2}}: Fix some $\rho\in (0,1)$ from Lemma \ref{L:5.2}. Thus the probability of $\{N_0+N_m> \rho n\}$ converges to zero uniformly with respect to $\theta$ and $F$. Therefore it suffices to show that the mean total variation distance between ${\cal L}_{\theta,F}^D(\Phi\mid N)$ and $P_{\theta,F}^C(\cdot\mid N)$, restricted to the event ${\cal N} := \{N_0+N_m \leq \rho n\}$, tends to zero uniformly in $\theta$ and $F$ as well. The first (conditional) probability measure has the Lebesgue density
\begin{align*} h_{\theta}(x\mid N) \, = \, (2\pi)^{1/2-1/m} (&\det \Delta \Psi_N(\vartheta))^{1/2} (\det \Delta \Psi_N(x))^{-1} \\
& \cdot \exp\big\{-\big((\nabla \Psi_N)^{-1}(-x) -\vartheta\big)^T \Delta \Psi_N(\vartheta) \big((\nabla \Psi_N)^{-1}(-x) - \vartheta\big) / 2\big\}\,, \end{align*}
on the range ${\cal R}$ of $-\nabla \Psi_N$, on which $h_\theta(\cdot\mid N)$ is supported and on which the function $\nabla \Psi_N$ has an inverse mapping. We write $g_\theta(\cdot\mid N)$ for the density of $P_{\theta,F}^C(\cdot \mid N) = \mathsf{N}(-\nabla \Psi_N(\vartheta),\Delta \Psi_N(\vartheta))$. Moreover note that 
$$ E \, 1_{\cal N} \, \int |h_\theta(x\mid N) - g_\theta(x\mid N)| dx \, \leq \, 2 E\, 1_{\cal N} \, \int_{\cal R} \big|h_\theta(x\mid N) - g_\theta(x\mid N)\big| dx\,. $$
Applying the integral substitution via $-\nabla \Psi_N$ the right hand side of the above inequality equals $2\, E 1_{\cal N} |Y-1|$ where
\begin{align} \nonumber Y &\, := \, \frac{\det \Delta \Psi_N(X)}{\det \Delta \Psi_N(\vartheta)} \, \exp\Big\{ -\frac12 \Big( \big(\nabla \Psi_N(X) - \nabla \Psi_N(\vartheta)\big)^T \{\Delta \Psi_N(\vartheta)\}^{-1}\\  \nonumber &  \hspace{5.5cm}  \big(\nabla \Psi_N(X) - \nabla \Psi_N(\vartheta)\big) - (X-\vartheta)^T \Delta \Psi_N(\vartheta) (X-\vartheta)\Big)\Big\}\,, \end{align}
where ${\cal L}(X\mid N) = \mathsf{N}(\vartheta,\{\Delta \Psi_N(\vartheta)\}^{-1})$. All third-order partial derivatives of $\Psi_N$ are bounded by $6n$ so that
$$ \nabla \Psi_N(X) - \nabla \Psi_N(\vartheta) \, = \, \Delta \Psi_N(\vartheta)\, (X-\vartheta) \, + \, R_1\,, $$
where the remainder term $R_1$ satisfies $|R_1| \leq 6n m^{3/2} |X-\vartheta|^2$. The matrix-valued function $\Delta \Psi_N$ has the following Lipschitz property (with respect to the Frobenius norm),
$$ \big\| \Delta \Psi_N(X) - \Delta \Psi_N(\vartheta) \big\|_F \, \leq \, 6 n m^{3/2} |X-\vartheta|\,. $$ 
The Theorem of Courant-Fischer yields that
$$ \sup_{j=1,\ldots,m-1} \big|\lambda_j(X) - \lambda_j(\vartheta)\big| \, \leq \, 6 n m^{3/2} |X-\vartheta|\,, $$ 
where $\lambda_j(X)$ and $\lambda_j(\vartheta)$ denote the eigenvalues of the matrices $\Delta \Psi_N(X)$ and $\Delta \Psi_N(\vartheta)$, respectively, in decreasing order. We learn from the Lemmata \ref{L:4.3} and \ref{L:5.1} that 
$$ \inf_n\, (m_n^2/n)\cdot \inf_{\theta} \inf_{j} \lambda_j(\vartheta) \, > \, 0\,, $$
for $m=m_n$. Thus, on the event ${\cal C} := {\cal N} \cap \{|X-\vartheta| \leq \alpha_n m_n^{-9/2}\}$, for any sequence $(\alpha_n)\downarrow 0$, we deduce that
$$ \Big|1 - \frac{\det \Delta \Psi_N(X)}{\det \Delta \Psi_N(\vartheta)}\Big| \, = \, \Big|1 - \prod_{j=1}^{m_n-1}\Big(1 + \frac{\lambda_j(X) - \lambda_j(\vartheta)}{\lambda_j(\vartheta)}\Big)\Big|\, \leq \, \mbox{const.}\cdot \alpha_n\,, $$
where, in the sequel, const. stands for a constant only depending on $\rho$ and $R$. Furthermore,
\begin{align*} & \big(\nabla \Psi_N(X) - \nabla \Psi_N(\vartheta)\big)^T \{\Delta \Psi_N(\vartheta)\}^{-1} \big(\nabla \Psi_N(X) - \nabla \Psi_N(\vartheta)\big) - (X-\vartheta)^T \Delta \Psi_N(\vartheta) (X-\vartheta) \\
& \, = \, 2 R_1^T (X-\vartheta) + R_1^T \{\Delta \Psi_N(\vartheta)\}^{-1} R_1 \, \leq \, \mbox{const.}\cdot \big( n \alpha_n^3 m_n^{-12} + n \alpha_n^4 m_n^{-13}\big)\,,
\end{align*}
 holds true on the event ${\cal C}$. Any selection of $(\alpha_n)_n$ such that $$\lim_{n\to\infty} \alpha_n m_n^{-4} n^{1/3} = 0\,, $$ guarantees uniform convergence of $E \, 1_{\cal C} \, |Y-1|$ to zero. On the other hand the probability of ${\cal N} \backslash \{|X-\vartheta| \leq \alpha_n m_n^{-9/2}\}$ also tends to zero uniformly with respect to $\theta$ and $F$ if 
$$ \alpha_n n^{1/2} m_n^6 \, \to \, \infty\,, $$
as $n\to\infty$ since ${\cal L}(X-\vartheta) = \mathsf{N}(0,\{\Delta \Psi_N(\vartheta)\}^{-1})$. As $\sup_n m_n^\beta n < \infty$ for some $\beta>13$ such a choice of $(\alpha_n)_n$ exists. Then, 
$$ \lim_{n\to\infty} \sup_{\theta,F} E 1_{\cal N} |Y \cdot 1_{\cal C} - 1| \, = \, 0\,. $$
As $Y$ is non-negative, $E 1_{\cal N} Y \leq 1$ and $\lim_{n\to\infty} \sup_{\theta,F} (1-P_{\theta,F}^C({\cal N})) = 0$ we arrive at 
$$ \lim_{n\to\infty} \sup_{\theta,F} E 1_{\cal N} |Y - 1| \, = \, 0\,, $$
which completes the proof. \hfill $\square$ \\

\noindent {\it Proof of Lemma \ref{L:5.3}}: Setting
$$ \eta_k(\beta) \, := \, \exp\{\beta-\theta_k\} / \big(1 + \exp\{\beta-\theta_k\}\big)\,, \qquad k=1,\ldots,m\,, $$
we may write
$$ q_k(\theta,F) \, = \, \int \big\{*_{k=1}^m \mathsf{B}(1,\eta_k(\beta))\big\}\, dF(\beta)\,. $$
As $\theta \in [-R,R]^m$ we have that
\begin{align*} \eta_k(\beta) & \, \leq \, \exp(2R) \cdot \eta_1(\beta)\,, \\
1 -\eta_k(\beta) & \, \leq \, \exp(2R) \cdot\big(1 - \eta_1(\beta)\big)\,,
\end{align*}
for all $\beta\in \mathbb{R}$ and $k=2,\ldots,m$. 

On the sets $A_1 := \{\beta \, : \, \eta_1(\beta) \leq c m^{-3/4}\}$ and $A_2 := \{\beta \, : \, 1 - \eta_1(\beta) \leq c m^{-3/4}\}$ for some constant $c>0$, we apply Poisson approximation of binomial distributions. Precisely, an inequality of Le Cam (see p. 657 in DasGupta (2008), for instance) yields that 
$$ \mathsf{TV}\big(\mathsf{P}(Q(\beta)),*_{k=1}^m \mathsf{B}(1,\eta_k(\beta))\big) \, \leq \, 2 \exp(4R) c^2 m^{-1/2}\,, \qquad \forall \beta\in A_1\,,$$
where $Q(\beta) \, := \sum_{k=1}^m \eta_k(\beta)$. Put $I_k(\delta) := Q^{-1}([k-\delta,k+\delta])$ for $k=1,\ldots,m-1$; $I_0(\delta) := Q^{-1}([0,\delta])$ and $I_m(\delta) := Q^{-1}([m-\delta,m-\delta/2])$ for some fixed $\delta\in (0,1)$. By Stirling's approximation, 
\begin{align*} b(m,k,\beta) & \, := \, \big\{*_{l=1}^m \mathsf{B}(1,\eta_l(\beta))\big\}(k) \, \geq \, \exp\{-Q(\beta)\} Q^k(\beta) / k! - 2 \exp(4R) \cdot c^2 m^{-1/2}  \\ & 
\, \geq \, \exp\big\{-1/(12k) - \delta\big\} (1-\delta/k)^k / \sqrt{2\pi k} - 2 \exp(4R) \cdot c^2 m^{-1/2} \\ &
\, \geq \, \mbox{const.}\cdot m^{-1/2}\,, \end{align*}
for all $\beta \in A_1 \cap I_k(\delta)$, $k\geq 1$, and a constant factor only depending on $R$, when choosing the constant $c>0$ sufficiently small. For $k=0$ this bound applies as well. 

For $\beta \in A_2 \cap I_k(\delta)$ the identical lower bound applies since $q_k(\theta,F)$ is viewed as the density of $*_{l=1}^m \mathsf{B}(1,1-\eta_l(\beta))$ at $m-k$. 

Finally we consider the complement $A_3 := \mathbb{R} \backslash (A_1\cup A_2)$. Clearly $\eta_1(\beta) \in (c m^{-3/4}, 1 - cm^{-3/4})$ holds for all $\beta \in A_3$. By Fourier inversion,
\begin{align*}
b(m,k,\beta) & \, = \, \frac1{2\pi} \int_{-\pi}^\pi \exp\big\{-it (k-Q(x))\big\} f_m(t,\beta) dt\,,
\end{align*}
with \begin{align*}
f_m(t,\beta) & \, := \, \prod_{l=1}^m \exp\{i t (1-\eta_l(\beta))\} \cdot \eta_l(\beta) + \exp\{-i t \eta_l(\beta))\} \cdot (1 - \eta_l(\beta)) \,, \\
\eta(\beta) & \, := \, \sum_{l=1}^m \eta_l(\beta) (1-\eta_l(\beta))\,. 
\end{align*}
As in the proof of Lemma \ref{L:4.1} we derive that
$$ |f_m(t,\beta)| \, \leq \, \exp\Big\{-\frac2{\pi^2} \eta(\beta)\Big\}\,, $$
for all $t\in [-\pi,\pi]$. Moreover, for all $\beta \in A_3$, we have
$$ \eta(\beta) \, \geq \, c \exp\{-2R\} m^{1/4} (1 - c \exp\{-2R\} m^{-3/4})\,. $$
Put $\nu := \{D (\log \eta(\beta))/\eta(\beta)\}^{1/2}$ for some constant $D>0$ sufficiently large. Then $|f_m(t,\beta)| \leq m^{-D/(2\pi^2)}$ if $|t| \in (\nu,\pi]$. Otherwise, for $t\in [-\nu,\nu]$, Taylor approximation yields that $$ f_m(t,\beta) = \exp\big\{- t^2 \eta(\beta) / 2\big\} \cdot \big(1 + \Delta_m(t,\beta)\big)\,, $$
with the remainder $|\Delta_m(t,\beta)| \leq \mbox{const.}\cdot \nu$ for some universal constant factor. Combining these facts with
$$ \big|1 - \exp\big\{-it (k-Q(\beta))\big\}\big| \, \leq \, |t| \delta\,, $$
for all $\beta \in I_k(\delta)$, we deduce that 
\begin{equation} \label{eq:eta} b(m,k,\beta) \, \geq \, \mbox{const.}\cdot m^{-1/2}\,, \end{equation}
for all $\beta \in I_k(\delta)\cap A_3$ and some universal positive constant. Summarisingly, the inequality (\ref{eq:eta}) has been verified for all $\beta\in I_k(\delta)$ where the constant factor only depends on $R$ and $\delta$. 

We conclude that
$$ q_k(\theta,F) \, \geq \, \mbox{const.}\cdot m^{-1/2} \cdot \int_{I_k(\delta)} f(x) dx\,, $$
for all $k=0,\ldots,m$ where the constant does not depend on $k$. As the derivative of $Q$ is bounded by $m$ the length of the interval $I_k(\delta)$ has the lower bound $\delta/(2m)$. Moreover,
$$ \sup\{|x| : x\in I_k(\delta)\} \, \leq \, R + |\log(\delta/2)| + \log m\,, $$
for any $\delta \in (0,1/2)$ and all $k=0,\ldots,m$. Finally, by the tail condition (\ref{eq:tail}) on $f$, the proof is completed. \hfill $\square$ \\

\noindent {\it Proof of Lemma \ref{L:5.4}}: The total variation distance between $\mathsf{N}(n \tilde{q}(\theta,F), n \tilde{Q}(\theta,F) - n \tilde{q}(\theta,F) \tilde{q}(\theta,F)^T)$ and ${\cal L}(\tau(N^{**})\mid V)$ is bounded from above by
\begin{align*}
&1_{(-\infty,\zeta)}(V) \, + \, 1_{[\zeta,\infty)}(V) \\ & \cdot \mathsf{TV}\big(\mathsf{N}\big(n \tilde{q}(\theta,F), n^3 (\tilde{Q}(\theta,F) - \tilde{q}(\theta,F) \tilde{q}(\theta,F)^T) / V^2\big),\mathsf{N}\big(n \tilde{q}(\theta,F), n (\tilde{Q}(\theta,F) - \tilde{q}(\theta,F) \tilde{q}(\theta,F)^T)\big)\big) \\ & \, \leq \, 1_{(-\infty,\zeta)}(V)\, + \, 1_{[\zeta,\infty)}(V) \cdot \mathsf{H}\big(\mathsf{N}\big(0, (n^2/V^2) I\big) , N(0,I)\big) \\ &
 \, \leq \, 1_{(-\infty,\zeta)}(V)\, + \, \sqrt{2(m+1)} \cdot 1_{[\zeta,\infty)}(V) \cdot \big|n^2/V^2 - 1\big|\,,
\end{align*}
where $\mathsf{H}$ denotes the Hellinger distance and $I$ the $(m+1)\times (m+1)$-identity matrix. Therein equation (A.4) in Rei{\ss} (2011) has been used to bound the Hellinger distance between normal distributions. Applying the expectation to the above term we obtain
\begin{align*}
P^F_{\theta,F}(V < \zeta) & \, + \, \sqrt{2(m+1)} \zeta^{-2} E \big|n^2 - V^2\big| \, \leq \, 4/n \, + \, 4 \sqrt{2(m+1)} (2 n^{3/2} + n) / n^2\,, 
\end{align*}
as ${\cal L}_{\theta,F}^F(V) \sim N(n,n)$ and $\zeta = n/2$. Thanks to the conditions on $m_n$ in Theorem \ref{T:5}, the above expression tends to zero uniformly in $\theta\in\Theta$ and $F\in {\cal F}$. \hfill $\square$ \\

\noindent {\it Proof of Theorem \ref{T:appl1}}: We consider that
\begin{align} \nonumber
& P^D_{\theta,F}\big(\theta \not\in \hat{E} \mid N\big) \\ \nonumber 
& \, \leq 1_{[0,n)}(N_0+N_m)  P^D_{\theta,F}\big(\varepsilon^T \Delta \Psi_{N}(\vartheta)^{-1/2} \Delta \Psi_{N}(T^{**})  \Delta \Psi_{N}(\vartheta)^{-1/2} \varepsilon > \iota \mid N\big)  +  P_{\theta,F}^D(N_0+N_m = n) \\ \nonumber
& \, \leq 1_{[0,n)}(N_0+N_m) P^D_{\theta,F}\big(|\varepsilon|^2 \cdot\big\{1 + 6 n m^{3/2} \big\|\{\Delta \Psi_{N}(\vartheta)\}^{-1/2}\big\|^3 |\varepsilon|\big\} > \iota  \mid N\big)   +  P_{\theta,F}^D(N_0+N_m = n),
\end{align}
where $\|\cdot\|$ denotes the usual matrix norm which is induced by the Euclidean metric; we have used that
\begin{equation*} \sup_{x \in \mathbb{R}^{m-1}}\, \sup_{l,j,j'} \Big| \frac{\partial^3}{\partial x_l \partial x_j \partial x_{j'}} \Psi_{N}(x) \Big| \, \leq \, 6n\,. \end{equation*}
By the Lemmata \ref{L:4.3} and \ref{L:5.1} we deduce that
\begin{equation*} P^D_{\theta,F}\big(\theta \not\in \hat{E}\big)  \leq 2 P^D_{\theta,F}\big(N_0 + N_m > \rho n\big) + P\big(|\varepsilon| > h_n n^{1/2} m^{-9/2}\big) + P\big(\big|\varepsilon\big|^2 \cdot\big(1 + 6 c (1-\rho)^{-3/2} h_n\big) > \iota \big)\,, \end{equation*}
for some sequence $(h_n)_n$ such that $(h_n m_n)_n \downarrow 0$ (with $m=m_n$) but $\big(h_n^2 n m_n^{-10}\big)_n \uparrow \infty$, some constant $c>0$ only depending on $R$ and some fixed $\rho\in (0,1)$ from Lemma \ref{L:5.2}. Taking the supremum over $\theta\in\Theta$ and $F\in {\cal F}$ and, then, the limit superior $n\to\infty$ on both sides of the above inequality we conclude that
$$ \limsup_{n\to\infty} \sup_{\theta,F} P^D_{\theta,F}\big(\theta \not\in \hat{E}\big) \, \leq \, \limsup_{n\to\infty} P\big(\big|\varepsilon\big|^2 > \iota \big) \, = \, 1-\alpha\,, $$
where we have used Lemma \ref{L:5.2} and the fact that the $\chi^2(m-1)$-density, as a consecutive sequence of convolutions, is bounded uniformly with respect to $m$ where $\iota \asymp m$ (since $\alpha\in (0,1)$ is fixed). Moreover the choice of $\iota$ is crucial in the last step. 

The maximal axis $\hat{e}$ of $\hat{E}$ turns out to be $2/\sqrt{\lambda_{\min}}$ where $\lambda_{\min}$ is the smallest positive eigenvalue of the matrix $Z (Z^\dagger Z)^{-1} \Delta \Psi_{N}(T^{**}) (Z^\dagger Z)^{-1} Z^\dagger$. For all $x \in \mathbb{R}^m$ with $\sum_{j=1}^m x_j = 0$ we have $\big|(Z^\dagger Z)^{-1} Z^\dagger x\big| \geq |x|$. Therefore $\lambda_{\min}$ is bounded from below by the smallest eigenvalue $\lambda'_{\min}$ of the matrix $\Delta \Psi_{N}(T^{**})$. Then it follows from the Lemmata \ref{L:4.3} and \ref{L:5.1} that
$$ \lambda'_{\min} \, \geq \, \mbox{const.}\cdot \big\{(1-\rho) n/m_n^2 \, - \, (1-\rho)^{-1/2} |\varepsilon| \cdot m_n / n^{1/2}\big\}\,, $$
  holds on the event $\{N_0+N_m \leq \rho n\}$. Using Lemma \ref{L:5.2} we establish that
$$ 1/\lambda_{\min} \, = \, {\cal O}_P\big(m_n^2/n\big)\,, $$
uniformly with respect to $\theta$ and $F$. \hfill $\square$

\end{document}